\documentclass{amsart}

\usepackage{amsfonts,amssymb,amscd,amstext}

\usepackage{graphicx}

\renewcommand{\leq}{\leqslant}
\renewcommand{\geq}{\geqslant}

\renewcommand{\ge}{\geqslant}
\newcommand{\ptl}{\partial}

\newcommand{\rr}{{\mathbb{R}}}
\newcommand{\la}{\lambda}

\newcommand{\hh}{{\mathbb{H}}}

\newcommand{\sph}{{\mathbb{S}}}
\newcommand{\pp}{\mathcal{P}}

\newcommand{\sub}{\subset}
\newcommand{\subeq}{\subseteq}

\newcommand{\escpr}[1]{\big<#1\big>}
\newcommand{\Sg}{\Sigma}
\newcommand{\Om}{\Omega}
\newcommand{\eps}{\varepsilon}
\newcommand{\var}{\varphi}
\newcommand{\ga}{\gamma}
\newcommand{\Ga}{\Gamma}

\newcommand{\nuh}{\nu_{H}}

\newcommand{\nh}{N_{H}}

\DeclareMathOperator{\divv}{div}

\newtheorem{theorem}{Theorem}[section]
\newtheorem{proposition}[theorem]{Proposition}
\newtheorem{lemma}[theorem]{Lemma}
\newtheorem{corollary}[theorem]{Corollary}

\theoremstyle{definition}

\newtheorem{remark}[theorem]{Remark}
\newtheorem{example}[theorem]{Example}

\theoremstyle{remark}

\newenvironment{enum}{\begin{enumerate}
}{\end{enumerate}}

\numberwithin{equation}{section}

\begin{document}

\title[Area-stationary surfaces in the Heisenberg
group]{Area-stationary surfaces \\ in the Heisenberg group $\hh^1$}

\author[M.~Ritor\'e]{Manuel Ritor\'e} \address{Departamento de
Geometr\'{\i}a y Topolog\'{\i}a \\
Universidad de Granada \\ E--18071 Granada \\ Espa\~na}
\email{ritore@ugr.es}

\author[C.~Rosales]{C\'esar Rosales} \address{Departamento de
Geometr\'{\i}a y Topolog\'{\i}a \\
Universidad de Granada \\ E--18071 Granada \\ Espa\~na}
\email{crosales@ugr.es}

\date{\today}

\thanks{Both authors have been supported by MCyT-Feder research
project MTM2004-01387} 

\begin{abstract}
We use variational arguments to introduce a notion of mean curvature
for surfaces in the Heisenberg group $\hh^1$ endowed with its
Carnot-Carath\'eodory distance.  By analyzing the first variation of
area, we characterize $C^2$ area-stationary surfaces as those
with mean curvature zero (or constant if a volume-preserving condition
is assumed) and such that the characteristic curves meet orthogonally
the singular curves.  Moreover, a Minkowski type formula
relating the area, the mean curvature, and the volume is obtained for
volume-preserving area-stationary surfaces enclosing a given region.

As a consequence of the characterization of area-stationary surfaces,
we refine the Bernstein type theorem given in \cite{chmy} and
\cite{gp} to describe entire area-stationary graphs over the
$xy$-plane in $\hh^1$.  A calibration argument shows that these
graphs are globally area-minimizing.

Finally, by using the description of the singular set in \cite{chmy},
the characterization of area-stationary surfaces, and the ruling
property of constant mean curvature surfaces, we prove our main
results where we classify volume-preserving area-stationary surfaces
in $\hh^1$ with non-empty singular set.  In particular, we deduce the
following counterpart to Alexandrov uniqueness theorem in Euclidean
space: any compact, connected, $C^2$ surface in $\hh^1$,
area-stationary under a volume constraint, must be congruent with a
rotationally symmetric sphere obtained as the union of all the
geodesics of the same curvature joining two points.  As a consequence,
we solve the isoperimetric problem in $\hh^1$ assuming $C^2$
smoothness of the solutions.
\end{abstract}

\maketitle

\thispagestyle{empty}

\section{\textbf{Introduction}}
\label{sec:intro}
\setcounter{equation}{0}

In the last years the study of variational questions in sub-Riemannian
geometry has received an increasing interest.  In particular, the
desire to achieve a better understanding of global variational
questions involving the area, such as the \emph{Plateau problem} or
the \emph{isoperimetric problem}, has motivated the recent development
of a theory of \emph{constant mean curvature surfaces} in the
\emph{Heisenberg group} $\hh^1$ endowed with its
\emph{Carnot-Carath\'eodory distance}.

It is well-known that constant mean curvature surfaces arise as
critical points of the area for variations preserving the volume
enclosed by the surface.  In this paper, we are interested in surfaces
immersed in the Heisenberg group which are \emph{stationary points} of
the sub-Riemannian area, with or without a \emph{volume constraint}.
In order to precise the situation and state our results we recall some
facts about the Heisenberg group, that will be treated in more detail
in Section~\ref{sec:preliminaries}.

We denote by $\hh^1$ the $3$-dimensional \emph{Heisenberg group},
which we identify with the Lie group $\mathbb{C}\times\rr$, where the
product is given by
\[
[z,t]*[z',t']=[z+z',t+t'+\text{Im}\big(z\overline{z}'\big)].
\]
The Lie algebra of $\hh^1$ is generated by three left invariant vector
fields $\{X,Y,T\}$ with one non-trivial bracket relation given by
$[X,Y]=-2T$.  The $2$-dimensional distribution generated by $\{X,Y\}$
is called the \emph{horizontal distribution} in $\hh^1$.  Usually
$\hh^1$ is endowed with a structure of sub-Riemannian manifold by
considering the Riemannian metric on the horizontal distribution so
that the basis $\{X,Y\}$ is orthonormal.  This metric allows to
measure the length of horizontal curves and to define the
\emph{Carnot-Carath\'eodory distance} between two points as the
infimum of length of horizontal curves joining both points, see
\cite{Gr}.  It is known that the Carnot-Carath\'edory distance can be
approximated by the distance functions associated to a family of
dilated Riemannian metrics, see \cite{gromov2}, \cite{pansu} and
\cite[\S 1.10]{mont}.  The Heisenberg group $\hh^1$ is also a
pseudo-hermitian manifold.  It is the simplest one and can be seen as
a blow-up of general pseudo-hermitian manifolds
(\cite[Appendix]{chmy}).  In addition, $\hh^1$ is also a \emph{Carnot
group} since its Lie algebra is stratified and $2$-nilpotent, see
\cite{dgn}.

Since $\hh^1$ is a group one can consider
its Haar measure, which turns out to coincide with the Lebesgue
measure in $\rr^3$.  From the notions of distance and measure one can
also define the Minkowski content and the sub-Riemannian perimeter of
a set, and the spherical Hausdorff measure of a surface, so that
different surface measures may be given on $\hh^1$.  As it is shown in
\cite{msc} and \cite{fsc}, all these notions of ``area'' coincide for
a $C^2$ surface.

In this paper we introduce the notions of volume and area in $\hh^1$
as follows.  We consider the left invariant Riemannian metric
$g=\escpr{\cdot\,,\cdot}$ on $\hh^1$ so that $\{X,Y,T\}$ is an
orthonormal basis at every point.  We define the volume $V(\Om)$ of a
Borel set $\Om\subeq\hh^1$ as the Riemannian measure of the set.  The
area of an immersed $C^1$ surface $\Sg$ in $\hh^1$ is defined as the
integral
\[
A(\Sg)=\int_{\Sg}|N_{H}|\,d\Sg,
\]
where $N$ is a unit vector normal to the surface, $N_{H}$ denotes the
orthogonal projection onto the horizontal distribution, and $d\Sg$ is
the Riemannian area element induced on $\Sg$ by the metric $g$.  This
definition of area agrees for $C^2$ surfaces with the ones mentioned
above.

With these notions of volume and area, we study in
Section~\ref{sec:meancurvature} surfaces in $\hh^1$ which are
\emph{stationary} points of the area either for arbitrary variations,
or for variations preserving the volume enclosed by the surface.  As
in Riemannian geometry, one may expect that some geometric quantity
defined on such a surface vanishes or remains constant.  By using the
first variation of area in Lemma~\ref{lem:dp/dt} we will see that any
$C^2$ area-stationary surface under a volume constraint must have
\emph{constant mean curvature}.  The mean curvature $H$ of a surface
$\Sg$ is defined in \eqref{eq:mc} as the Riemannian divergence
relative to $\Sg$ of the \emph{horizontal unit normal vector} to $\Sg$
given by $\nuh=N_{H}/\,|N_{H}|$.  We remark that a notion of mean
curvature in $\hh^1$ for graphs over the $xy$-plane was previously
introduced by S.~Pauls \cite{pauls}.  A more general definition of
mean curvature has been proposed by J.-H.~Cheng, J.-F.~Hwang,
A.~Malchiodi and P.~Yang \cite{chmy}, and by N.~Garofalo and S.~Pauls
\cite{gp}.  As was shown in \cite{revolucion} our definition agrees
with all the previous ones.

The analysis of the \emph{singular set} plays an important role in the
study of area-stationary surfaces in $\hh^1$.  Given a surface $\Sg$
immersed in $\hh^1$, the singular set $\Sg_0$ of $\Sg$ is the set of
points where $\Sg$ is tangent to the horizontal distribution.  Its
structure has been determined for surfaces with bounded
mean curvature in \cite{chmy}, where it is proved that $\Sg_{0}$
consists of isolated points and $C^1$ curves, see Theorem
\ref{th:chmy} for a more detailed description.  The \emph{regular
part} $\Sg-\Sg_{0}$ of $\Sg$ is foliated by horizontal curves called
the \emph{characteristic curves}.  As is pointed out in \cite{chmy},
when the surface $\Sg$ has constant mean curvature $H$, any of
these curves is part of a \emph{geodesic} in $\hh^1$ of curvature $H$.
In particular, any surface in $\hh^1$ with $H\equiv 0$ is foliated, up
to the singular set, by horizontal straight lines.

The recent study of surfaces with constant mean curvature in $\hh^1$
has mainly focused on minimal surfaces (those with $H\equiv 0$).  In
fact, many interesting questions of the classical theory of minimal
surfaces in $\rr^3$, such as the Plateau problem, the Bernstein
problem, or the global behavior of properly embedded surfaces, have
been treated in $\hh^1$, see \cite{pauls}, \cite{chmy}, \cite{gp},
\cite{ch}, and \cite{pauls2}.  These works also provide a rich variety
of examples of minimal surfaces.  However, in spite of the last
advances, very little is known about non-minimal constant mean curvature
surfaces in $\hh^1$.  It is easy to check that a graph $t=u(x,y)$ of
class $C^2$ in $\hh^1$ with constant mean curvature $H$ satisfies the
following degenerate (elliptic and hyperbolic) PDE
\[
(u_{y}+x)^2u_{xx}-2\,(u_{y}+x)(u_{x}-y)\,u_{xy}+(u_{x}-y)^2u_{yy}=
-2H\,((u_{x}-y)^2+(u_{y}+x)^2)^{3/2}.
\]
In \cite{chmy} some relevant properties concerning the above equation, such
as the uniqueness of solutions for the Dirichlet problem or the
structure of the singular set, are studied.  As to the examples, the
only known complete surfaces with non-vanishing constant mean
curvature are the compact spherical ones described in \cite{monti} and
\cite{leomas}, and the complete surfaces of revolution that we
classified in \cite{revolucion}.

Now we briefly describe the organization and the results obtained in
this paper.  After the preliminaries Section~\ref{sec:preliminaries},
we make a detailed study of sub-Riemannian geodesics and Jacobi fields
in Section~\ref{sec:geodesics}.  In Section~\ref{sec:meancurvature} we
look at the first variation of area and prove a Minkowski-type formula
for an area-stationa\-ry surface under a volume constraint relating
area, volume and the mean curvature, Theorem~\ref{th:minkowski}.
Then, a detailed analysis of the first variation of area, together
with the aforementioned description of the singular set in
Theorem~\ref{th:chmy}, leads us to prove in Theorem \ref{th:constant}
that an immersed surface is area-stationary if and only if its mean
curvature is zero (or constant under a volume constraint) and the
characteristic curves meet orthogonally the singular curves.  This
result allows us to refine in Section \ref{sec:bernstein} the
Bernstein-type theorems given in \cite{chmy} and \cite{gp} for minimal
graphs in $\hh^1$.  We classify all entire area-stationary graphs in
$\hh^1$ over the $xy$-plane in Theorem~\ref{th:bernstein}, and show that
they are globally area-minimizing in Theorem~\ref{th:areaminimizing}.
In Section \ref{sec:mainresult}, we prove our main results, where we
completely describe immersed area-stationary surfaces in $\hh^1$ under
a volume constraint with non-vanishing mean curvature and non-empty
singular set, Theorems~\ref{th:spheres} and \ref{th:classification}.
As a consequence we deduce an Alexandrov uniqueness type theorem for
compact surfaces, Theorem~\ref{th:alexandrov}, and we solve the
isoperimetric problem in $\hh^1$ assuming $C^2$ regularity of the
solutions in Theorem~\ref{th:iso}.

Now we describe our results in more detail.

A classical formula by Minkowski in Euclidean space involving the
integral of the support function over a compact surface in $\rr^3$
with constant mean curvature yields the relation $A=3HV$, where $A$
is the area of the surface, $V$ is the volume enclosed, and $H$ is the
mean curvature of the surface.  Our analysis of the first variation of
the sub-Riemannian area and the existence in $\hh^1$ of a
one-parameter group of dilations provide a Minkowski-type formula for
a surface $\Sg$ which is area-stationary under a volume constraint in
$\hh^1$.  Such a formula reads
\[
3A=8HV,
\]
where $A$ is the sub-Riemannian area of $\Sg$, $H$ the mean curvature
of $\Sg$, and $V$ the volume enclosed.

From previous works, as \cite{chmy}, \cite{dgn}, \cite{gp}, and
\cite{revolucion}, it was already known that a necessary condition for
a surface $\Sg$ to be area-stationary is that the mean curvature of
$\Sg$ must be zero (or constant if the surface is area-stationary
under a volume constraint).  In Theorem \ref{th:constant} we show that
such a condition is not sufficient.  To obtain a stationary point for
the area we must require in addition that the \emph{characteristic
curves meet orthogonally the singular curves}.  We prove this result
by obtaining an expression for the first variation of area for
arbitrary variations of the surface $\Sg$, not only for those fixing
the singular set.  Observe that the situation is different from the one
in Riemannian geometry, where stationary surfaces are precisely those
with vanishing mean curvature.

As a consequence of this analysis, we show that most of the entire
graphs obtained in \cite{chmy} and \cite{gp} with mean curvature zero
are not area-stationary.  We refine their result to prove that
the only entire area-stationary graphs over the $xy$-plane in $\hh^1$
are the Euclidean planes and vertical rotations of the graphs
\[
u(x,y)=xy+(ay+b),
\]
where $a$, $b\in\rr$.  Geometrically, the latter surfaces can be
described as the union of all horizontal lines in $\hh^1$ which are
orthogonal to a given horizontal line (the singular curve).  By using
a calibration argument, we can prove that they are globally
area-minimizing.  This result is similar to the Euclidean one, where
planes, the only entire minimal graphs in $\rr^3$, are
area-minimizing.  In \cite[\S 6]{chmy}, also by a calibration
argument, it was proved that a compact portion of the regular part of
a graph with mean curvature zero is area-minimizing.

It was already known that the regular part of a surface $\Sg$ immersed
in $\hh^1$ with constant mean curvature $H$ is foliated by horizontal
geodesics of curvature $H$.  We derive in Section~\ref{sec:geodesics}
an intrinsic equation for such geodesics and for Jacobi fields, and
show in Theorem \ref{th:ruled} that the characteristic curves of the
surface are geodesics of curvature $H$.  This is the starting point,
together with the local description of the singular set in
Theorem~\ref{th:chmy}, to construct new examples and to classify
surfaces of constant mean curvature in $\hh^1$.

In Section \ref{sec:mainresult} we use this idea to describe any
complete, volume-preserving area-stationary surface $\Sg$ in $\hh^1$
with non-vanishing mean curvature and non-empty singular set.  We
prove in Theorem~\ref{th:spheres} that if $\Sg$ has at least one
isolated singular point then it must be congruent with one of the
compact spherical examples $\sph_\la$ obtained as the union of all the
geodesics of curvature $\la>0$ joining two given points
(Example~\ref{ex:spheres}).  Then, we introduce in
Proposition~\ref{prop:sigmala} a procedure to construct examples of
complete surfaces with non-vanishing constant mean curvature $\la$.
Geometrically these surfaces consist of a horizontal curve $\Ga$ in
$\hh^1$, from which geodesics of curvature $\la$ leave (or enter)
orthogonally.  An analysis of the variational vector field associated
to this family of geodesics is necessary to understand the behavior of
the geodesics far away from $\Ga$.  It follows that the resulting
surface has two singular curves apart from $\Ga$.  Moreover, the
family of geodesics meets both curves orthogonally if and only if they
are equidistant to $\Ga$.  This geometric property allows to conclude
in Theorem~\ref{th:curve} the strong restriction that \emph{the
singular curves of any volume-preserving area-stationary surface in
$\hh^1$ with $H\neq 0$ are geodesics of $\hh^1$}.  This is the key
ingredient to classify in Theorem~\ref{th:classification} all surfaces
of this kind.  It follows that they must be congruent either with the
cylindrical embedded surfaces in Example~\ref{ex:cilindros} or with
the helicoidal immersed surfaces in Example~\ref{ex:helices}.

This technique can also be used to describe complete area-stationary
surfaces with singularities.  It was proved in
\cite[Proposition~2.1]{ch} and \cite[Lemma~8.2]{gp} that Euclidean
planes are the only complete minimal surfaces in $\hh^1$ with at least
one isolated singular point.  In Proposition~\ref{prop:onesingu} we
give a nice geometric description of complete area-stationary surfaces
with singular curves: the singular curve is a unique, arbitrary
horizontal curve and the surface consists of the union of all the
horizontal lines orthogonal to this singular curve.

Alexandrov uniqueness theorem in Euclidean space states that the only
embedded compact surfaces with constant mean curvature in $\rr^3$ are
round spheres.  This result is not true for immersed surfaces as
illustrated by the toroidal examples in \cite{wente}.  In
pseudo-hermitian geometry, an interesting restriction on the topology
of an immersed compact surface with bounded mean curvature inside a
$3$-spherical pseudo-hermitian manifold was given in \cite{chmy},
where it was proved that such a surface is homeomorphic either to a
sphere or to a torus.  As shown in \cite{chmy} this bound on the genus
is optimal on the standard pseudo-hermitian $3$-sphere, where examples
of constant mean curvature spheres and tori may be given.  This
estimate on the genus is also valid in $\hh^1$ since the proof is
based on the local description of the singular set
(Theorem~\ref{th:chmy}) and on the Hopf Index Theorem.  In
Theorem~\ref{th:alexandrov} we prove the following counterpart in
$\hh^1$ to Alexandrov uniqueness theorem in $\rr^3$: any compact,
connected, $C^2$ immersed volume-preserving area-stationary surface
$\Sg$ in $\hh^1$ is congruent with a sphere $\sph_\la$.  In particular
we deduce the non-existence of volume-preserving area-stationary
immersed tori in $\hh^1$.

Finally in Section~\ref{sec:iso} we study the \emph{isoperimetric
problem} in $\hh^1$.  This problem consists of finding sets in $\hh^1$
minimizing the sub-Riemannian perimeter under a volume constraint.  It
was proved by G. P. Leonardi and S. Rigot \cite{lr} that the solutions
to this problem exist and they are bounded, connected open sets.  This
information is clearly far from characterizing isoperimetric sets.  In
the last years many authors have tried to adapt to the Heisenberg
group different proofs of the classical isoperimetric inequality in
Euclidean space.  In \cite{monti2}, \cite{monti} and \cite{leomas} it
was shown that there is no a direct counterpart in $\hh^1$ to the
Brunn-Minkowski inequality in Euclidean space, with the surprising
consequence that the Carnot-Carath\'edory balls in $\hh^1$ cannot be
the solutions.  Recently, interest has focused on solving the
isoperimetric problem restricted to certain sets with additional
symmetries.  It was proved by D. Danielli, N. Garofalo and D.-M. Nhieu
that the sets bounded by the spherical surfaces $\sph_\la$ are the
unique solutions in the class of sets bounded by two $C^1$ radial
graphs over the $xy$-plane \cite[Theorem~14.6]{dgn}.  In
\cite{revolucion} we pointed out that assuming $C^2$ smoothness and
rotationally symmetry of isoperimetric regions, these must be
congruent with the spheres $\sph_\la$.  We finish this work by showing
in Theorem~\ref{th:iso} that the spherical surfaces $\sph_\la$ are the
unique isoperimetric regions in $\hh^1$ assuming $C^2$ regularity of
the solutions, solving a conjecture by P.~Pansu \cite[p.~172]{pansu3}.
Regularity of solutions is still a hard, open question.

After the distribution of this paper, we have noticed some related
works.  In \cite{chy}, interesting results for graphs in the
Heisenberg group $\hh^n$ have been established.  In particular, the
authors prove in \cite[p.~30]{chy} that $C^2$ minimal graphs in
$\hh^1$ are area-minimizing if and only if the characteristic curves
meet orthogonally the singular curves.  In \cite{dgn2} it is proved
that the sets bounded by the spheres $\sph_{\la}$ are the unique
isoperimetric regions in the class of sets bounded by the union of two
$C^1$ graphs over the $xy$-plane.  In \cite{dgn3} the authors show
that there exists a family of entire \emph{intrinsic minimal graphs}
in $\hh^1$ that are not area-minimizing.  In \cite{bc} the mean
curvature flow of a $C^2$ convex surface in $\hh^1$, described as the
union of two radial graphs, is proved to converge to a sphere
$\sph_{\la}$.  In \cite{bscv}, it is introduced a general calibration
method to study the Bernstein problem for entire regular intrinsic
minimal graphs in the Heisenberg group $\hh^n$.  Finally we mention
the interesting survey \cite{survey}, where the authors give a broad
overview of the isoperimetric problem in $\hh^n$.

\section{{\bf Preliminaries}}
\label{sec:preliminaries}
\setcounter{equation}{0}

The \emph{Heisenberg group} $\hh^1$ is the Lie group $(\rr^3,*)$,
where the product $*$ is defined, for any pair of points $[z,t]$,
$[z',t']\in\rr^3\equiv\mathbb{C}\times\rr$, as
\[
[z,t]*[z',t']:=[z+z',t+t'+\text{Im}(z\overline{z}')], \qquad (z=x+iy).
\]
For $p\in\hh^1$, the \emph{left translation} by $p$ is the
diffeomorphism $L_p(q)=p*q$.  A basis of left invariant vector fields
(i.e., invariant by any left translation) is given by
\begin{equation*}
X:=\frac{\ptl}{\ptl x}+y\,\frac{\ptl}{\ptl t}, \qquad
Y:=\frac{\ptl}{\ptl y}-x\,\frac{\ptl}{\ptl t}, \qquad
T:=\frac{\ptl}{\ptl t}.
\end{equation*}
The \emph{horizontal distribution} $\mathcal{H}$ in $\hh^1$ is the
smooth planar one generated by $X$ and $Y$.  The \emph{horizontal
projection} of a vector $U$ onto $\mathcal{H}$ will be denoted by
$U_{H}$.  A vector field $U$ is called \emph{horizontal} if $U=U_H$.
A \emph{horizontal curve} is a $C^1$ curve whose tangent vector lies
in the horizontal distribution.

We denote by $[U,V]$ the Lie bracket of two $C^1$ vector fields $U$, $V$
on $\hh^1$. Note that $[X,T]=[Y,T]=0$, while
$[X,Y]=-2T$.  The last equality implies that $\mathcal{H}$ is a
bracket generating distribution. Moreover, by Frobenius Theorem we
have that $\mathcal{H}$ is nonintegrable.  The vector fields $X$ and $Y$
generate the kernel of the (contact) $1$-form
$\omega:=-y\,dx+x\,dy+dt$.

We shall consider on $\hh^1$ the (left invariant) Riemannian metric
$g=\escpr{\cdot\,,\cdot}$ so that $\{X,Y,T\}$ is an orthonormal basis
at every point, and the associated Levi-Civit\'a connection $D$.  The
modulus of a vector field $U$ will be denoted by $|U|$.  The following
derivatives can be easily computed
\begin{alignat}{2}
\notag D_{X}X&=0, \qquad \ \ \ \, D_{Y}Y=0, \qquad \,D_{T}T=0, \\
\label{eq:christoffel}
D_{X}Y&=-T, \qquad \, D_{X}T=Y, \qquad \, D_{Y}T=-X, \\
\notag D_{Y}X&=T, \qquad \ \ \,\,D_{T}X=Y, \qquad D_{T}Y=-X.
\end{alignat}
For any vector field $U$ on $\hh^1$ we define $J(U):=D_UT$.  Then we
have $J(X)=Y$, $J(Y)=-X$ and $J(T)=0$, so that $J^2=-\text{Identity}$
when restricted to the horizontal distribution.  It is also clear that
\begin{equation}
\label{eq:conmute}
\escpr{J(U),V}+\escpr{U,J(V)}=0,
\end{equation}
for any pair of vector fields $U$ and $V$.  The endomorphism $J$
restricted to the horizontal distribution is an involution of
$\mathcal{H}$ that, together with the contact $1$-form
$\omega=-y\,dx+x\,dy+dt$, provides a pseudo-hermitian structure on
$\hh^1$, as stated in the Appendix in \cite{chmy}.

Let $\ga:I\to\hh^1$ be a piecewise $C^1$ curve defined on a compact
interval $I\sub\rr$.  The \emph{length} of $\ga$ is the usual
Riemannian length $L(\ga):=\int_{I}|\dot{\ga}|$, where $\dot{\ga}$ is
the tangent vector of $\ga$.  For two given points in $\hh^1$ we can
find, by Chow's connectivity Theorem \cite[p.  95] {Gr}, a horizontal
curve joining these points.  The \emph{Carnot-Carath\'edory distance}
$d_{cc}$ between two points in $\hh^1$ is defined as the infimum of
the length of horizontal curves joining the given points.

Now we introduce notions of volume and area in $\hh^1$.  The volume
$V(\Om)$ of a Borel set $\Om\subeq\hh^1$ is the Riemannian volume of
the left invariant metric $g$, which coincides with the Lebesgue
measure in $\rr^3$.  Given a $C^1$ surface $\Sg$ immersed in $\hh^1$,
and a unit vector field $N$ normal to $\Sg$, we define the area of
$\Sg$ by
\begin{equation}
\label{eq:area}
A(\Sg):=\int_{\Sg}|N_{H}|\,d\Sg,
\end{equation}
where $N_{H}=N-\escpr{N,T}\,T$, and $d\Sg$ is the Riemannian area
element on $\Sg$.  If $\Sg$ is a $C^1$ surface enclosing a bounded set
$\Om$ then $A(\Sg)$ coincides with the $\hh^1$-perimeter of $\Om$, as
defined in \cite{cng}, \cite{fsc} and \cite{revolucion}.  The area of
$\Sg$ also coincides with the Minkowski content in $(\hh^1,d_{cc})$ of a set
$\Om\subset\hh^1$ bounded by a $C^2$ surface $\Sg$, as proved in
\cite[Theorem~5.1]{msc}, and with the 3-dimensional spherical Hausdorff
measure in $(\hh^1,d_{cc})$ of $\Sg$, see \cite[Corollary~7.7]{fsc}.

For a $C^1$ surface $\Sg\sub\hh^1$ the \emph{singular set} $\Sg_0$
consists of those points $p\in\Sg$ for which the tangent plane
$T_p\Sg$ coincides with the horizontal distribution.  As $\Sg_0$ is
closed and has empty interior in $\Sg$, the \emph{regular set}
$\Sg-\Sg_0$ of $\Sg$ is open and dense in $\Sg$.  It was proved in
\cite[Lemme 1]{d2}, see also \cite[Theorem~1.2]{balogh}, that the
Hausdorff dimension with respect to the Riemannian distance on
$\hh^1$ of $\Sg_{0}$ is less than two.

If $\Sg$ is a $C^1$ oriented surface with unit normal vector $N$, then
we can describe the singular set $\Sg_0\sub\Sg$, in terms of $N_H$, as
$\Sg_{0}=\{p\in\Sg:N_H(p)=0\}$.  In the regular part $\Sg-\Sg_0$, we
can define the \emph{horizontal unit normal vector} $\nu_H$, as in
\cite{dgn}, \cite{revolucion} and \cite{gp} by
\begin{equation}
\label{eq:nuh}
\nu_H:=\frac{N_H}{|N_H|}.
\end{equation}
Consider the \emph{characteristic vector field} $Z$ on $\Sg-\Sg_0$
given by
\begin{equation}
\label{eq:zeta}
Z:=J(\nu_H).
\end{equation}
As $Z$ is horizontal and orthogonal to $\nu_H$, we conclude that $Z$
is tangent to $\Sg$.  Hence $Z_{p}$ generates the intersection of
$T_{p}\Sg$ with the horizontal distribution.  The integral curves of
$Z$ in $\Sg-\Sg_0$ will be called \emph{characteristic curves} of
$\Sg$.  They are both tangent to $\Sg$ and horizontal. Note that these
curves depend on the unit normal $N$ to $\Sg$.   If we define
\begin{equation}
\label{eq:ese}
S:=\escpr{N,T}\,\nu_H-|N_H|\,T,
\end{equation}
then $\{Z_{p},S_{p}\}$ is an orthonormal basis of $T_p\Sg$ whenever
$p\in\Sg-\Sg_0$.

In the Heisenberg group $\hh^1$ there is a one-parameter group of
\emph{dilations} $\{\varphi_s\}_{s\in\rr}$ generated by the vector
field
\begin{equation}
\label{eq:w}
W:=xX+yY+2tT.
\end{equation}
From the Christoffel symbols \eqref{eq:christoffel}, it can be easily
proved that $\divv W=4$, where $\divv W$ is the Riemannian divergence
of the vector field $W$.  We may compute $\varphi_{s}$ in coordinates
to obtain
\begin{equation}
\label{eq:dilations}
\varphi_{s}(x_{0},y_{0},t_{0})=(e^sx_{0},e^sy_{0},e^{2s}t_{0}).
\end{equation}
From this expression we get, for fixed $s$ and $p\in\hh^1$, that
$(d\varphi_{s})_{p}(X_{p})=e^s X_{\varphi_{s}(p)}$,
$(d\varphi_{s})_{p}(Y_{p})=e^s Y_{\varphi_{s}(p)}$, and
$(d\varphi_{s})_{p}(T_{p})=e^{2s} T_{\varphi_{s}(p)}$.

Any isometry of $(\hh^1,g)$ leaving invariant the horizontal
distribution preserves the area of surfaces in $\hh^1$.  Examples of
such isometries are left translations, which act transitively on
$\hh^1$.  The Euclidean rotation of angle $\theta$ about the $t$-axis
given by
\[
(x,y,t)\mapsto r_{\theta}(x,y,t)=(\cos\theta\,x-\sin\theta\,y,
\sin\theta\,x+\cos\theta\,y,t),
\]
is also an area-preserving isometry in $(\hh^1,g)$ since it transforms the
orthonormal basis $\{X,Y,T\}$ at the point $p$ into the orthonormal
basis $\{\cos\theta\,X+\sin\theta\,Y, -\sin\theta\,X+\cos\theta\,Y,
T\}$ at the point $r_{\theta}(p)$.

\section{\textbf{Geodesics and Jacobi fields in the Heisenberg group $\hh^1$}}
\label{sec:geodesics}
\setcounter{equation}{0}

 Usually, geodesics in $\hh^1$ are defined as horizontal curves whose
 length coincides with the Carnot-Carath\'eodory distance between its
 endpoints.  It is known that geodesics in $\hh^1$ are curves of class
 $C^\infty$, see \cite[Lemma~2.5]{monti2}.  We are interested in
 computing the equations of geodesics in terms of geometric data of
 the left invariant metric $g$ in $\hh^1$.  For that we shall think of
 a geodesic in $\hh^1$ as a smooth horizontal curve that is a critical
 point of length under any variation by horizontal curves with fixed
 endpoints.  In this section we will obtain an \emph{intrinsic}
 equation for the geodesics in terms of the left invariant metric $g$.

Let $\ga:I\to\hh^1$ be a $C^2$ horizontal curve defined on a compact
interval $I\sub\rr$. A variation of $\ga$ is a $C^2$ map
$F:I\times J\to\hh^1$, where $J$ is an open interval around the
origin, such that $F(s,0)=\ga(s)$.  We denote $\ga_\eps(s)=F(s,\eps)$.
Let $V_{\eps}(s)$ be the vector field along $\ga_{\eps}$ given by $(\ptl
F/\ptl\eps)(s,\eps)$.  Trivially $[V_{\eps},\dot{\ga_{\eps}}]=0$. Let
$V=V_{0}$.  We say that
the variation is \emph{admissible} if the curves $\ga_\eps$ are
horizontal and have fixed boundary points.  For such a variation it is
clear that $V$ vanishes at the endpoints of $\ga$.  Moreover, we have
$\escpr{\dot{\ga_{\eps}},T}=0$.  As a consequence
\begin{align*}
0=\frac{d}{d\eps}\bigg|_{\eps=0}\escpr{\dot{\ga_{\eps}},T}
&=\escpr{D_{V}\dot{\ga_{\eps}},T}+\escpr{\dot{\ga},D_{V}T}
\\
&=\escpr{D_{\dot{\ga}}V,T}+\escpr{\dot{\ga},J(V)}
\\
&=\dot{\ga}\big(\escpr{V,T}\big)-\escpr{V,D_{\dot{\ga}}T}
+\escpr{\dot{\ga},J(V_H)}
\\
&=\dot{\ga}\big(\escpr{V,T}\big)-\escpr{V_{H},J(\dot{\ga})}
+\escpr{\dot{\ga},J(V_{H})}
\\
&=\dot{\ga}\big(\escpr{V,T}\big)-2\,\escpr{V_{H},J(\dot{\ga})},
\end{align*}
where in the last equality we have used \eqref{eq:conmute}.

Conversely, if $V$ is a $C^1$ vector field along $\ga$ vanishing at the
endpoints and satisfying the equation
\begin{equation}
\label{eq:condition}
\dot{\ga}\big(\escpr{V,T}\big)=2\,\escpr{V_{H},J(\dot{\ga})},
\end{equation}
then it is easy to check that there is an admissible variation of
$\ga$ so that the associated vector field coincides with $V$.  Indeed,
since $V=f\dot{\ga}+V_{0}$, with $V_{0}\perp \dot{\ga}$, we may assume
that $V$ is orthogonal to $\ga$.  Define, for $s\in I$ and $\eps$
small, $F(s,\eps):=\exp_{\ga(s)}(\eps\,V(s))$, where $\exp$ is the
exponential map associated to the Riemannian metric $g$ in $\hh^1$.
If $V$ is horizontal in some interval of $\ga$ then, by
\eqref{eq:condition}, we have $V=V_{H}=\la\dot{\ga}$, so that $V$
vanishes.  If $V(s_{0})$ is not horizontal, $F$ defines locally a
surface which is transversal to the horizontal distribution.  This
surface is foliated by horizontal curves.  So there is a $C^2$ function
$f(s,\eps)$ such that $\ga_{\eps}(s):=\exp_{\ga(s)}(f(s,\eps)\,V(s))$
is a horizontal curve.  We may take $f$ so that $(\ptl
f/\ptl\eps)(s_{0},0)=1$.  The vector field $V_{1}$ associated to the
variation by horizontal curves $\ga_{\eps}$, is given by $(\ptl
f/\ptl\eps)(s,0)\,V(s)$, and satisfies equation
\eqref{eq:condition}.  Since $V$ also satisfies this equation we
obtain that $(\ptl^2f/\ptl s\,\ptl\eps)(s,0)=0$, and $(\ptl
f/\ptl\eps)(s,0)$ is constant.  As $(\ptl f/\ptl\eps)(s_{0},0)=1$ we
conclude that $V_{1}(s)=V(s)$.

\begin{proposition}
Let $\ga:I\to\hh^1$ be a $C^2$ horizontal curve parameterized by
arc-length.  Then $\ga$ is a critical point of length for any
admissible variation if and only if there is $\lambda\in\rr$ such that
$\ga$ satisfies the second order ordinary differential equation
\begin{equation}
\label{eq:geodesic}
D_{\dot{\ga}}\dot{\ga}+2\lambda\,J(\dot{\ga})=0.
\end{equation}
\end{proposition}

\begin{proof}
Let $V$ be the vector field of an admissible variation $\gamma_\eps$
of $\gamma$.  Since $\ga$ is parameterized by arc-length, by the first
variation of length \cite[\S1,(1.3)]{cheeger}, we know that
\begin{equation}
\label{eq:lprima}
\frac{d}{d\eps}\bigg|_{\eps=0} L(\ga_{\eps})=-
\int_{I}\escpr{D_{\dot{\ga}}\dot{\ga},V}.
\end{equation}
Suppose that $\gamma$ is a critical point of length for any admissible
variation.  As $|\dot{\ga}|=1$ we deduce that
$\escpr{D_{\dot{\ga}}\dot{\ga},\dot{\ga}}=0$.  On the other hand, as
$\ga$ is a horizontal curve, we have
$\escpr{D_{\dot{\ga}}\dot{\ga},T}=0$.  So $D_{\dot{\ga}}\dot{\ga}$ is
proportional to $J(\dot{\ga})$ at any point of $\ga$.  Assume, without
loss of generality, that $I=[0,a]$.  Consider a $C^1$ function $f:I\to\rr$
vanishing at the endpoints and such that $\int_{I}f=0$.  Let $V$ be
the vector field on $\ga$ so that $V_{H}=f\,J(\dot{\ga})$ and
$\escpr{V,T}(s)=2\,\int_{0}^s f$.  As $V$ satisfies
\eqref{eq:condition}, inserting it in the first variation of length
\eqref{eq:lprima}, we obtain
\[
\int_{I}f\,\escpr{D_{\dot{\ga}}\dot{\ga},J(\dot{\ga})}=0.
\]
As $f$ is an arbitrary $C^1$ mean zero function we conclude that
$\escpr{D_{\dot{\ga}}\dot{\ga},J(\dot{\ga})}$ is constant.  Hence we find
$\lambda\in\rr$ so that $\ga$ satisfies equation \eqref{eq:geodesic}.  The
proof of the converse is easy taking into account \eqref{eq:lprima} and
\eqref{eq:condition}.
\end{proof}

We will say that a $C^2$ horizontal curve $\ga$ is a \emph{geodesic of
curvature} $\la$ if it is parameterized by arc-length and satisfies
equation~\eqref{eq:geodesic}.  Observe that the parameter $\la$ in
\eqref{eq:geodesic} changes to $-\la$ for the reversed curve
$\ga(-t)$.

Given a point $p\in\hh^1$, a unit horizontal vector $v\in T_{p}\hh^1$, and
$\lambda\in\rr$, we denote by $\ga_{p,v}^\lambda$ the unique
solution to \eqref{eq:geodesic} with initial conditions $\ga(0)=p$,
$\dot{\ga}(0)=v$. Note that $\ga_{p,v}^\lambda$ is a geodesic since it
is horizontal and parameterized by arc-length (the functions
$\escpr{\dot{\ga},T}$ and $|\dot{\ga}|^2$ are constant along any
solution of \eqref{eq:geodesic}).

Let us now compute the equation of the geodesics in
Euclidean coordinates.  Consider a $C^2$ curve $\ga(s)=(x(s),y(s),t(s))$
parameterized by arc-length.  Then
\[
\dot{\ga}=(\dot{x},\dot{y},\dot{t})=
\dot{x}\,X+\dot{y}\,Y+(-\dot{x}y+x\dot{y}+\dot{t})\,T,
\]
so that $\ga$ is horizontal if and only if
\[
-\dot{x}y+x\dot{y}+\dot{t}=0.
\]
Moreover:
\[
D_{\dot{\ga}}\dot{\ga}=\ddot{x}\,X+\ddot{y}\,Y, \qquad
2\lambda\,J(\dot{\ga})=2\lambda\,\big(\dot{x}\,Y-\dot{y}\,X\big).
\]
Hence $\ga=(x,y,t)$ is a geodesic of curvature $\lambda$ if it
satisfies the following system of equations
\begin{align*}
\ddot{x}&=2\lambda\,\dot{y}, \\
\ddot{y}&=-2\lambda\,\dot{x}, \\
\dot{t}&=\dot{x}y-x\dot{y}.
\end{align*}
\indent Let us solve first the case $\lambda\neq 0$.  Calling
$\dot{x}=u$, $\dot{y}=v$ we get $\ddot{u}+(2\lambda)^2u=0$, from
which, if $(\dot{x}(0),\dot{y}(0))=(A,B)$, we have $u(0)=A$,
$\dot{u}(0)=2\lambda B$, and
\begin{align*}
\dot{x}(s)&=u(s)=A\,\cos(2\lambda\,s)+B\,\sin(2\lambda\,s), \\
\dot{y}(s)&=v(s)=-A\,\sin(2\lambda\,s)+B\,\cos(2\lambda\,s).
\end{align*}
If $(x(0),y(0),t(0))=(x_{0},y_{0},t_{0})$, then:
\begin{align}
\nonumber
x(s)&=x_{0}+A\,\bigg(\frac{\sin(2\lambda\,s)}{2\lambda}\bigg)
+B\,\bigg(\frac{1-\cos(2\lambda\,s)}{2\lambda}\bigg),
\\
\label{eq:geocoor}
y(s)&=y_{0}-A\,\bigg(\frac{1-\cos(2\lambda\,s)}{2\lambda}\bigg)
+B\,\bigg(\frac{\sin(2\lambda\,s)}{2\lambda}\bigg),
\\
\nonumber
t(s)&=t_{0}+\frac{1}{2\lambda}\,
\bigg(s-\frac{\sin(2\lambda\,s)}{2\lambda}\bigg)
\\
\nonumber &+(Ax_{0}+By_{0})\left(\frac{1-\cos(2\lambda\,
s)}{2\la}\right)
-(Bx_{0}-Ay_{0})\,\left(\frac{\sin(2\lambda\, s)}{2\la}\right),
\end{align}
which are Euclidean helices of vertical axis.  Thus, we have recovered
the expressions in \cite[p.~28]{andre} and \cite[p.~160]{monti2}.
Assume now that $\lambda=0$.  In this case, we have the following
system of ordinary differential equations
\begin{align*}
\ddot{x}&=0, \\
\ddot{y}&=0, \\
\dot{t}&=\dot{x}y-x\dot{y}.
\end{align*}
For initial conditions $(x(0),y(0),t(0))=(x_{0},y_{0},t_{0})$,
$\dot{x}(0)=A$, $\dot{y}(0)=B$, we get
\begin{align*}
x(s)&=x_{0}+As, \\
y(s)&=y_{0}+Bs, \\
t(s)&=t_{0}+(Ay_{0}-Bx_{0})\,s,
\end{align*}
which are Euclidean horizontal lines. This fact was previously observed in
\cite[Proposition~4.1]{chmy}. We conclude that complete geodesics in
$\hh^1$ are horizontal lifts of curves with constant geodesic
curvature in the Euclidean $xy$-plane (circles or straight lines).

\begin{remark}
\label{re:propertiesgeo}
1.  Any isometry in $(\hh^1,g)$ preserving the horizontal distribution
transforms geodesics in geodesics since it respects the Levi-Civit\'a
connection and commutes with $J$.

2. A dilation $\varphi_{s}(x,y,t)=(e^sx,e^sy,e^{2s}t)$ carries
geodesics of curvature $\la$ to geodesics of curvature $e^{-s}\la$.

3.  If we consider the geodesic $\ga_{0,v}^\lambda$, where $v$ is a
horizontal unit vector in $T_{0}\hh^1$ and $\lambda\neq 0$, then the
coordinate $t(s)$ in \eqref{eq:geocoor} is monotone increasing
and unbounded.  It follows that $\ga_{0,v}^\lambda$ leaves every
compact set in finite time.  The same is true for any other horizontal
geodesic, since it can be transformed into $\ga_{0,v}^\la$ by a left
translation.
\end{remark}

\begin{lemma}
Let $\lambda>0$, $p\in\hh^1$, and $v$, $w\in T_{p}\hh^1$ horizontal unit
vectors with $v\neq w$.  Then
$\ga_{p,v}^\lambda(\pi/\lambda)=\ga_{p,w}^\lambda(\pi/\lambda)$ and
$\ga_{p,v}^\lambda(s_{1})\neq\ga_{p,w}^\lambda(s_{2})$ for all $s_{1}$,
$s_{2}\in (0,\pi/\lambda)$.
\end{lemma}

\begin{proof}
After applying a left translation and a rotation about the $t$-axis we
may assume that $p=(0,0,0)$, that $v=(1,0,0)$ and that
$w=(\cos\theta,\sin\theta,0)$, with $\cos\theta\neq 1$. From
\eqref{eq:geocoor}, we have that $\ga_{p,v}^\la$ is given by
\begin{align*}
x_{v}(s)&=(2\lambda)^{-1}\sin(2\lambda s), \\
y_{v}(s)&=(2\lambda)^{-1}\,\big(-1+\cos(2\lambda s)\big), \\
t_{v}(s)&=(2\lambda)^{-1}\,\big(s-(2\lambda)^{-1}\sin(2\lambda
s)\big),
\end{align*}
and $\ga_{p,w}^\lambda$ by
\begin{align*}
x_{w}(s)&=(2\lambda)^{-1}\,\big(\sin\theta+\sin(2\lambda s-\theta)\big), \\
y_{w}(s)&=(2\lambda)^{-1}\,\big(-\cos\theta+\cos(2\lambda s-\theta)\big), \\
t_{w}(s)&=(2\lambda)^{-1}\,\big(s-(2\lambda)^{-1}\sin(2\lambda
s)\big).
\end{align*}
Equality $\ga_{p,v}^\lambda(\pi/\lambda)=\ga_{p,w}^\lambda(\pi/\lambda)$ is
easily checked from these equations. Suppose that
$\ga_{p,v}^\lambda(s_{1})=\ga_{p,w}^\lambda(s_{2})$ for some
$s_1,s_2\in (0,\pi/\lambda)$. As $t_{v}=t_{w}$ is an increasing
function, we deduce $s_1=s_2$, and so there is $s\in (0,\pi/\la)$ such that
$(x_{v}(s),y_{v}(s))=(x_{w}(s),y_{w}(s))$. Therefore, we get
\begin{align*}
(1-\cos\theta)\,\sin(2\lambda s)&=(1-\cos(2\lambda s))\,\sin\theta,
\\
\sin(2\lambda s)\,\sin\theta&=(1-\cos\theta)\,(\cos(2\lambda s)-1),
\end{align*}
for some $s\in (0,\pi/\lambda)$. Finally, as the determinant
\[
\det\left(
\begin{array}{cc}
1-\cos\theta  &-\sin\theta
\\
\sin\theta  &1-\cos\theta
\end{array}
\right)\neq 0,
\]
we conclude that $\sin(2\lambda s)=0$ and $1-\cos(2 \lambda s)=0$, a
contradiction.
\end{proof}

\begin{figure}[h]
\centering{\includegraphics[width=6cm]{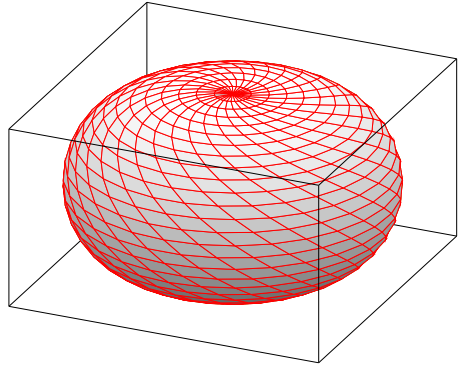}}
\caption{A spherical surface $\sph_{\la}$ given by the union of all
the geodesics of curvature $\la$ joining the poles.}
\label{fig:spheres}
\end{figure}

\begin{example}[Spheres in $\hh^1$]
\label{ex:spheres} Given $\lambda>0$, we define $\sph_{\lambda}$ as
the union of all geodesics $\ga_{0,v}^\lambda$ restricted to the
interval $[0,\pi/\lambda]$.  The lemma above implies that
$\sph_{\lambda}$ is a compact embedded surface homeomorphic to a
sphere, see Figure~\ref{fig:spheres}.  Any $\sph_{\la}$ has two
singular points at the \emph{poles} $(0,0,0)$ and $(0,0,\pi/(2\lambda^2))$.
Alternatively, it was proved in \cite[Proof of Theorem 3.3]{leomas}
that $\sph_{\la}$ can be described as the union of the following
radial graphs over the $xy$-plane
\begin{equation}
\label{eq:spheregraphs}
t=\frac{\pi}{2\la^2}\pm\frac{1}{2\la^2}\,\left (\la\rho\,\sqrt{1-\la^2\rho^2}
+\arccos (\la\rho)\right),\qquad\rho=\sqrt{x^2+y^2}\leq\frac{1}{\la}.
\end{equation}
From \eqref{eq:spheregraphs} we can see that $\sph_{\la}$ is $C^2$ but
not $C^3$ around the poles.  This was also observed in
\cite[Proposition~14.11]{dgn}.
\end{example}

Now, we prove some analytical properties for the vector field
associated to a variation of a curve which is a geodesic.

\begin{lemma}
\label{lem:jacobi1}
Let $\ga:I\to \hh^1$ be a geodesic of curvature $\la$, and $V$ the
$C^1$ vector field associated to a variation of $\ga$.  Then the function
\[
\la\,\escpr{V,T}+\escpr{V,\dot{\ga}}
\]
is constant along $\ga$.
\end{lemma}

\begin{proof}
First note that
\[
\dot{\ga}\,(\escpr{V,T})=\escpr{D_{\dot{\ga}}V,T}+\escpr{V,J(\dot{\ga})}
=\escpr{D_{V}\dot{\ga},T}-\escpr{\dot{\ga},J(V)}
=-2\,\escpr{\dot{\ga},J(V)},
\]
where we have used $[V,\dot{\ga}]=0$, equality \eqref{eq:conmute},
and that $\ga$ is a horizontal curve. On the other hand, we
have
\[
\dot{\ga}\,(\escpr{V,\dot{\ga}})=\escpr{D_{\dot{\ga}}V,\dot{\ga}}
+\escpr{V,-2\la\,J(\dot{\ga})}
=\escpr{D_{V}\dot{\ga},\dot{\ga}}+2\la\,\escpr{\dot{\ga},J(V)}
=2\la\,\escpr{\dot{\ga},J(V)},
\]
since $\ga$ is parameterized by arc-length and satisfies
\eqref{eq:geodesic}.  From the two equations above the result follows.
\end{proof}

As in Riemannian geometry we may expect that the vector field
associated to a variation of a given geodesic by geodesics of the same
curvature satisfies a certain second order differential equation.  In
fact, we have

\begin{lemma}
\label{lem:jacobi2}
Let $\ga_{\eps}$ be a variation of $\ga$ by geodesics of the same
curvature $\la$. Assume that the associated vector field $V$ is $C^2$.
Then $V$ satisfies
\begin{equation}
\label{eq:jacobieq}
\ddot{V}+\emph{R}(V,\dot{\ga})\dot{\ga}+
2\la\,(J(\dot{V})-\escpr{V,\dot{\ga}}\,T)=0,
\end{equation}
where $\emph{R}$ denotes the Riemannian curvature tensor in $(\hh^1,g)$.
\end{lemma}

\begin{proof}
As any $\ga_{\eps}$ is a geodesic of curvature $\la$, we have
\[
D_{\dot{\ga}_\eps}\dot{\ga_{\eps}}+2\la\,J(\dot{\ga_{\eps}})=0.
\]
Thus, if we derive with respect to $V$ and we take into account that
$D_VD_{\dot{\ga}}\dot{\ga}=\text{R}(V,\dot{\ga})\dot{\ga}+D_{\dot{\ga}}
D_V\dot{\ga}+D_{[V,\dot{\ga}]}\dot{\ga}$ and that $[V,\dot{\ga}]=0$, we deduce
\[
\ddot{V}+\text{R}(V,\dot{\ga})\dot{\ga}+2\la\,D_{V}J(\dot{\ga})=0.
\]
Finally, it is not difficult to see that
\[
D_{V}J(\dot{\ga})=J(D_{V}\dot{\ga})-\escpr{V,\dot{\ga}}\,T=
J(\dot{V})-\escpr{V,\dot{\ga}}\,T,
\]
and the proof follows.
\end{proof}

We call \eqref{eq:jacobieq} the \emph{Jacobi equation} for geodesics
in $\hh^1$ of curvature $\la$.  It is clearly a linear equation.  Any
solution of \eqref{eq:jacobieq} is a \emph{Jacobi field} along $\ga$.
It is easy to check that $V=f\dot{\ga}$ is a Jacobi field if and only
if $\ddot{f}\dot{\ga}+2\la \dot{f}J(\dot{\ga})=0$.  Thus, any tangent
Jacobi field to $\ga$ is of the form $(as+b)\,\dot{\ga}$, with $a=0$
when $\lambda\neq 0$.

\section{\textbf{Area-stationary surfaces. Minkowski formula in $\hh^1$}}
\label{sec:meancurvature}
\setcounter{equation}{0}

In this section we shall consider critical surfaces for the area
functional \eqref{eq:area} with or without a volume constraint.  Let
$\Sg$ be an oriented immersed surface of class $C^2$ in $\hh^1$.
Consider a $C^1$ vector field $U$ with compact support on $\Sg$.
Denote by $\Sg_{t}$, for $t$ small, the immersed surface
$\{\exp_{p}(tU_{p}); p\in\Sg\}$, where $\exp_{p}$ is the exponential
map of $(\hh^1,g)$ at the point $p$.  The family $\{\Sg_t\}$, for $t$
small, is the \emph{variation} of $\Sg$ induced by $U$. We remark
that our variations can move the singular set $\Sg_{0}$ of $\Sg$. Define
$A(t):=A(\Sg_t)$.  In case $\Sg$ is an embedded compact surface, it
encloses a region $\Om$ so that $\Sg=\ptl\Om$.  Let $\Om_{t}$ be the
region enclosed by $\Sg_{t}$ and define $V(t):=V(\Om_t)$.  We say that
the variation is \textit{volume-preserving} if $V(t)$ is constant for
$t$ small enough.  We say that $\Sg$ is \emph{area-stationary} if
$A'(0)=0$ for any variation of $\Sg$.  In case that $\Sg$ encloses a
bounded region, we say that $\Sg$ is \emph{area-stationary under a
volume constraint} or \emph{volume-preserving area-stationary} if
$A'(0)=0$ for any volume-preserving variation of $\Sg$.

Suppose that $\Om$ is the set bounded by a $C^2$ embedded compact
surface $\Sg=\ptl\Om$.  We shall always choose the unit \emph{inner}
normal $N$ to $\Sg$.  The computation of $V'(0)$ is well-known since
the volume is the one associated to a Riemannian metric, and we have
(\cite[\S 9]{simon})
\begin{equation}
\label{eq:1stvol}
V'(0)=\int_{\Om}\divv U\,dv=-\int_\Sg u\,d\Sg,
\end{equation}
where $u=\escpr{U,N}$, and $dv$ is the Riemannian volume element.  It
follows that $u$ has mean zero whenever the variation is
volume-preserving.  Conversely, it was proven in \cite[Lemma
2.2]{bdce} that, given a $C^1$ function $u:\Sg\to\rr$ with mean zero,
a volume-preserving variation of $\Om$ can be constructed so that the
normal component of the associated vector field equals $u$.

\begin{remark}
Let $\Sg$ be a $C^1$ compact immersed surface in $\hh^1$.  Observe
that the vector field $W$ defined in \eqref{eq:w} satisfies $\divv
W=4$, so that if $\Sg$ is embedded, the divergence theorem yields
\begin{equation}
\label{eq:immersedvolume}
\text{volume enclosed by
}\Sg=-\frac{1}{4}\int_{\Sg}\escpr{W,N}\,d\Sg,
\end{equation}
where $N$ is the inner unit normal to $\Sg$.  Formula
\eqref{eq:immersedvolume} can be taken as a definition for the volume
``enclosed'' by an oriented compact immersed surface in $\hh^1$.  The
first variation for this volume functional is given by
\eqref{eq:1stvol}.  Also the \emph{variation} of enclosed volume can
be defined for a noncompact surface.  We refer the reader to
\cite{bdce} for details.
\end{remark}

Now we will compute the first variation of area.  We need a previous
lemma.

\begin{lemma}
\label{lem:derivatives}
Let $\Sg\subset\hh^1$ be a $C^2$ surface and $N$ a unit vector normal
to $\Sg$.  Consider a point $p\in\Sg-\Sg_0$, the horizontal normal
$\nu_H$ defined in \eqref{eq:nuh}, and $Z=J(\nuh)$.  Then, for
any $u\in T_{p}\hh^1$ we have
\begin{align}
\label{eq:uno} D_{u}N_{H}&=(D_{u}N)_{H}-\escpr{N,T}\,J(u)
-\escpr{N,J(u)}\,T,
\\
\label{eq:dos} u\,(|N_H|)&=\escpr{D_{u}N, \nuh}-\escpr{N,T}\,
\escpr{J(u),\nuh},
\\
\label{eq:tres} D_{u}\nu_H&=|N_H|^{-1}\,
\big(\escpr{D_uN,Z}-\escpr{N,T}\,
\escpr{J(u),Z}\big)\,Z+\escpr{Z,u}\,T.
\end{align}
\end{lemma}

\begin{proof}
Equalities \eqref{eq:uno} and \eqref{eq:dos} are easily obtained since
$N_H=N-~\escpr{N,T}\,T$.  Let us prove \eqref{eq:tres}.  As $|\nu_H|=1$ and
$\{(\nuh)_{p},Z_{p},T_{p}\}$ is an orthonormal basis of $T_p\hh^1$, we get
\[
D_{u}\nuh=\escpr{D_{u}\nuh,Z}\,Z+ \escpr{D_{u}\nuh,T}\,T.
\]
Note that $\escpr{D_{u}\nuh,T}=-\escpr{\nuh,J(u)}=\escpr{Z,u}$ by
\eqref{eq:conmute}.  On the other hand, by using \eqref{eq:uno} and
the fact that $Z$ is tangent and horizontal, we deduce
\[
\escpr{D_{u}\nuh,Z}=|N_H|^{-1}\escpr{D_{u}N_H,Z}=
|N_H|^{-1}\big(\escpr{D_uN,Z}-\escpr{N,T}\escpr{J(u),Z}\big).\vspace{-0.75cm}
\]

\end{proof}

For a $C^1$ vector field $U$ defined on a surface $\Sg$, we denote by
$U^\top$ and $U^\bot$ the tangent and orthogonal projections,
respectively.  We shall also denote by $\divv_\Sg U$ the Riemannian
divergence of $U$ relative to $\Sg$, which is given by
$\divv_{\Sg}U(p):=\sum_{1=1}^{2}\escpr{D_{e_{i}}U,e_{i}}$ for any
orthonormal basis $\{e_{1},e_{2}\}$ of $T_{p}\Sg$.  Now, we can prove

\begin{lemma}
\label{lem:dp/dt}
Let $\Sg\subset\hh^1$ be an oriented $C^2$ immersed surface.  Suppose
that $U$ is a $C^1$ vector field with compact support on $\Sg$ and
normal component $u=\escpr{U,N}$.  Then the first derivative at $t=0$
of the area functional $A(t)$ associated to $U$ is given by
\begin{equation}
\label{eq:1stvar}
A'(0)=\int_\Sg
u\,\big(\divv_\Sg\nu_H\big)\,d\Sg
-\int_\Sg\divv_\Sg\big(u\,(\nu_H)^\top\big)\,d\Sg,
\end{equation}
provided $\divv_{\Sg}\nu_{H}\in L^1(\Sg)$.

Moreover, if $\Sg$ is area-stationary $($resp. volume-preserving
area-stationary$)$ then
\begin{equation}
\label{eq:1stvarst}
A'(0)=\int_\Sg u\,(\divv_\Sg\nu_H)\,d\Sg.
\end{equation}
\end{lemma}

\begin{proof}
First we remark that the Riemannian area of the singular set $\Sg_{0}$
of $\Sg$ vanishes, as was proved in \cite[Lemme 1]{d2} and
\cite[Theorem~1.2]{balogh}. Thus we can integrate over $\Sg$ functions
defined on the regular set $\Sg-\Sg_{0}$.

Let $\{\Sg_t\}$ be the variation of $\Sg$ associated to $U$, and let
$d\Sg_t$ be the Riemannian area element on $\Sg_t$.  Consider a $C^1$
vector field $N$ whose restriction to $\Sg_t$ coincides with a unit
vector normal to $\Sg_{t}$.  By using \eqref{eq:area} and the coarea
formula, we have
\[
A(t)=\int_{\Sg_t}|N_H|\,d\Sg_t=
\int_\Sg (|N_H|\circ\var_t)\,|\text{Jac}\,\var_t|\,d\Sg=
\int_{\Sg-\Sg_{0}}(|N_H|\circ\var_t)\,|\text{Jac}\,\var_t|\,d\Sg,
\]
where $\varphi_{t}(p)=\exp_{p}(tU_{p})$ and $\text{Jac}\,\var_t$ is
the Jacobian determinant of the map $\var_t:\Sg\to\Sg_{t}$.  Now, we
differentiate with respect to $t$, and we use the known fact that
$(d/dt)|_{t=0}\,|\text{Jac}\,\var_t|=\divv_\Sg U$
(\cite[\S 9]{simon}), to get
\begin{align*}
A'(0)&=\int_{\Sg-\Sg_{0}}\{U(|N_H|)+|N_H|\,\divv_{\Sg}U\}\,d\Sg
\\
&=\int_{\Sg-\Sg_{0}}\{U^\bot(|N_H|)+\divv_\Sg(|N_H|\,U)\}\,d\Sg
\\
&=\int_{\Sg-\Sg_{0}}\{\divv_{\Sg}(|N_H|\,U^\top)+
U^\bot(|N_H|)+|N_H|\,\divv_{\Sg}U^\bot\}\,d\Sg
\\
&=\int_{\Sg-\Sg_{0}}\{U^\bot(|N_H|)+|N_H|\,\divv_{\Sg}U^\bot\}\,d\Sg.
\end{align*}
To obtain the last equality we have used the Riemannian divergence
theorem to get that the integral of the divergence of the
Lipschitz vector field $|N_{H}|\,U^\top$ over $\Sg$ vanishes (the
modulus of a $C^1$ vector field in a Riemannian manifold is a
Lipschitz function). We observe that the function
$U^\bot(|N_H|)+|N_H|\,\divv_{\Sg}U^\bot$ is bounded in $\Sg-\Sg_{0}$
and so it lies in $L^1(\Sg)$.

On the other hand, we can use \eqref{eq:dos} to obtain
\[
U^{\bot}(|N_H|)=\escpr{D_{U^{\bot}}N,\nu_H} -\escpr{N,T}\,
\big<J(U^\bot),\nu_H\big>=-\escpr{\nabla_\Sg u,\nu_H},
\]
since $J(U^{\bot})$ is orthogonal to $\nu_H$ and
$D_{U^\bot}N=-\nabla_\Sg u$.  Here $\nabla_\Sg u$ represents the
gradient of $u$ relative to $\Sg$.  Then, we get in $\Sg-\Sg_{0}$
\begin{align*}
U^\bot(|N_H|)+|N_H|\,\divv_{\Sg}U^\bot
&=-(\nu_H)^{\top}(u)+u\,|N_H|\,\divv_{\Sg} N
\\
&=-\divv_\Sg\big(u\,(\nu_H)^\top\big)+u\,\divv_\Sg\big((\nu_H)^\top\big)
\\
&\quad+u\,\divv_\Sg(|N_H|\,N)
\\
&=-\divv_\Sg\big(u\,(\nu_H)^\top\big)+u\,\divv_\Sg\nu_H.
\end{align*}
As a consequence, we conclude that
\[
\int_{\Sg}\{U^\bot(|N_H|)+|N_H|\,\divv_{\Sg}U^\bot\}\,d\Sg=
\int_\Sg u\,\big(\divv_\Sg\nu_H\big)\,d\Sg-
\int_\Sg\divv_\Sg\big(u\,(\nu_H)^\top\big)\,d\Sg.
\]
Since we are assuming that $\divv_{\Sg}\nu_{H}\in L^1(\Sg)$ we
conclude that $\divv_{\Sg}(u\,(\nu_{H})^\top)\in L^1(\Sg)$ and so we
have
\[
A'(0)=\int_\Sg u\,\big(\divv_\Sg\nu_H\big)\,d\Sg-
\int_\Sg\divv_\Sg\big(u\,(\nu_H)^\top\big)\,d\Sg.
\]
Note that the second integral above vanishes by virtue of the
Riemannian divergence theorem whenever $u$ has compact support
disjoint from the singular set $\Sg_{0}$.

Now we shall prove \eqref{eq:1stvarst} for area-stationary surfaces
under a volume constraint.  The proof for area-stationary ones follows
with the obvious modifications.  Inserting in \eqref{eq:1stvar} mean
zero functions of class $C^1$ with compact support inside the regular set
$\Sg-\Sg_{0}$, we get that $\divv_{\Sg}\nuh$ is a constant function on
$\Sg-\Sg_{0}$.  If $u:\Sg\to\rr$ is any function, then we consider
$v:\Sg\to\rr$ with support in $\Sg-\Sg_{0}$ such that
$\int_{\Sg}(u+v)\,d\Sg=0$.  Inserting the mean zero function $u+v$ in
\eqref{eq:1stvar}, taking into account that $\divv_{\Sg}\nuh$ is
constant, and using the divergence theorem, we deduce that
$\int_{\Sg}\divv_{\Sg}(u\,(\nuh)^\top)\,d\Sg=0$, and
\eqref{eq:1stvarst} is proved.
\end{proof}

\begin{remark}
The first variation of area \eqref{eq:1stvarst} holds for any $C^2$
surface whenever the support of the vector field $U$ is disjoint from
the singular set, see also \cite[Lemma 3.2]{revolucion}. For
area-stationary surfaces we have shown that \eqref{eq:1stvarst} is
also valid for vector fields moving the singular set.
\end{remark}

For a $C^2$ immersed surface $\Sg$ in $\hh^1$ with a $C^1$ unit normal
vector $N$ we define, as in \cite{revolucion}, the \emph{mean
curvature} $H$ of $\Sg$ by the equality
\begin{equation}
\label{eq:mc}
-2H(p):=(\divv_{\Sg}\nuh)(p),\qquad p\in\Sg-\Sg_{0}.
\end{equation}
For any point in $\Sg-\Sg_{0}$ we consider the orthonormal basis of
the tangent space to $\Sg$ given by the vectors fields $Z$ and $S$
defined in \eqref{eq:zeta} and \eqref{eq:ese}.  Then we have
\[
-2H=\escpr{D_{Z}\nuh,Z}+\escpr{D_{S}\nu_{H},S}.
\]
From \eqref{eq:tres} in Lemma~\ref{lem:derivatives} we get
$\escpr{D_{S}\nuh,S}=0$, and we conclude that
\begin{equation}
\label{eq:dznuh}
-2H=\escpr{D_{Z}\nuh,Z}=|\nh|^{-1}\,\escpr{D_{Z}N,Z}.
\end{equation}

By using variations supported in the regular set of a surface immersed
in $\hh^1$, the first variation of area \eqref{eq:1stvar}, and the
first variation of volume \eqref{eq:1stvol}, we get

\begin{corollary}
\label{cor:hconstant}
Let $\Sg$ be a $C^2$ oriented immersed surface in $\hh^1$. Then
\begin{enum}
\item If $\Sg$ is area-stationary then the mean curvature of
$\Sg-\Sg_{0}$
vanishes.
\item If $\Sg$ is area-stationary under a volume constraint then the
mean curvature of $\Sg-\Sg_{0}$ is constant.
\end{enum}
\end{corollary}

\begin{remark}
The first derivative of area for variations with compact support in
the regular set, and the notion of mean curvature were given by
S.~Pauls~\cite{pauls} for graphs over the $xy$-plane in $\hh^1$, and
later extended by J.-H.~Cheng, J.-F.~Hwang, A.~Malchiodi and P.~Yang
\cite{chmy} for any surface inside a $3$-dimensional pseudo-hermitian
manifold.  The case of the $(2n+1)$-dimensional Heisenberg group
$\hh^n$ has been treated in \cite{dgn}, \cite{revolucion} and
\cite{bc}.  In \cite{hp}, R.~Hladky and S.~Pauls extend the notion of
mean curvature and Corollary \ref{cor:hconstant} for stationary
surfaces inside vertically rigid sub-Riemannian manifolds.  In the
recent paper \cite{chy} the first variation of area for graphs over
$\rr^{2n}$ has been computed for some more general variations moving
the singular set.  A definition of mean curvature by using Riemannian
approximations to the Carnot-Carath\'eodory distance in $\hh^1$ can be
found in \cite[p.  562]{ni} and \cite[\S 3]{survey}.
\end{remark}

\begin{example}
1. According to our definition, the graph of a $C^2$ function $u(x,y)$ has
constant mean curvature $H$ if and only if satisfies the equation
\[
(u_{y}+x)^2u_{xx}-2\,(u_{y}+x)(u_{x}-y)\,u_{xy}+(u_{x}-y)^2u_{yy}=-2H\,
((u_{x}-y)^2+(u_{y}+x)^2)^{3/2}
\]
outside the singular set.

2.  The spherical surface $\sph_{\la}$ in Example \ref{ex:spheres} has
constant mean curvature $\la$ with respect to the inner normal vector.
This can be seen by using the equation for constant mean curvature
graphs above and \eqref{eq:spheregraphs}.  It was proved in
\cite[Theorem~5.4]{revolucion} that $\sph_{\la}$ is, up to a vertical
translation, the unique $C^2$ compact surface of revolution around the
$t$-axis with constant mean curvature $\la$.
\end{example}

The ruling property of constant mean curvature surfaces in $\hh^1$,
already observed in \cite[(2.1), (2.24)]{chmy}, \cite[Corollary~5.3]{gp} and
\cite[Corollaries~4.5 and 6.10]{hp},
follows immediately from the expression \eqref{eq:dznuh} for the mean
curvature and the equation of geodesics \eqref{eq:geodesic}.

\begin{theorem}
\label{th:ruled}
Let $\Sg$ be an oriented immersed surface in $\hh^1$ of class $C^2$
with constant mean curvature $H$ outside the singular set.  Then any
characteristic curve of $\Sg$ coincides with an open arc of a
geodesic of curvature $H$.  As a consequence, the regular set of $\Sg$
is foliated by geodesics of curvature $H$.
\end{theorem}

\begin{proof}
A characteristic curve $\ga$ is parameterized by arc-length since the tangent to $\ga$ is the
characteristic vector field $Z$ defined in \eqref{eq:zeta}.  We must see that $\ga$ satisfies
equation \eqref{eq:geodesic} for $\lambda=H$. For any point of this curve, the vector fields
$Z$, $\nuh$ and $T$ provide an orthonormal basis of the tangent space to $\hh^1$.  Thus, we
have
\begin{align*}
D_{\dot{\ga}}\dot{\ga}=D_{Z}Z&=\escpr{D_{Z}Z,\nuh}\,\nuh+\escpr{D_{Z}Z,T}\,T
\\
&=-\escpr{Z,D_{Z}\nuh}\,\nuh-\escpr{Z,J(Z)}\,T
\\
&=2H\nuh=-2HJ(Z)=-2HJ(\dot{\ga}),
\end{align*}
where in the last equalities we have used \eqref{eq:dznuh} and that
$J(Z)=-\nuh$.
\end{proof}

\begin{remark}
Let $\Sg$ be a $C^2$ surface in $\hh^1$ and $\varphi_{s}$ the dilation of $\hh^1$ defined in
\eqref{eq:dilations}. The ruling property in Theorem~\ref{th:ruled} and the behavior of
geodesics under $\varphi_{s}$ (Remark~\ref{re:propertiesgeo}) imply that $\Sg$ has constant
mean curvature $\la$ if and only if the dilated surface $\varphi_{s}(\Sg)$ has constant mean
curvature $e^{-s}\la$.
\end{remark}

Now, we will prove a counterpart in $\hh^1$ of the Minkowski formula for compact surfaces in
$\rr^3$.  We need the following consequence of \eqref{eq:1stvarst}, Corollary
\ref{cor:hconstant} and the definition of the mean curvature

\begin{corollary}
\label{cor:caract} 
Let $\Sg\subset\hh^1$ be a $C^2$ surface enclosing a bounded region $\Om$.
Then $\Sg$ is volume-preserving area-stationary if and only if there is a real constant $H$
such that $\Sg$ is a critical point of the functional $A-2HV$ for any given variation.
\end{corollary}

This corollary and the existence in $\hh^1$ of a one-parameter group
of dilations allow us to prove the following Minkowski type formula
for volume-preserving stationary surfaces enclosing a bounded region
in $\hh^1$.  The result also holds for oriented compact immersed
surfaces in $\hh^1$ when the volume is given by
\eqref{eq:immersedvolume}.

\begin{theorem}[Minkowski formula in $\hh^1$]
\label{th:minkowski} Let $\Sg\subset\hh^1$ be a volume-preserving area-stationary $C^2$ surface
enclosing a bounded region $\Om$.  Then we have
\begin{equation}
\label{eq:minkowski}
3A(\Sg)=8H\,V(\Om),
\end{equation}
where $H$ is the mean curvature of $\Sg$ with respect to the inner normal vector.
\end{theorem}

\begin{proof}
We take the vector field $W$ in \eqref{eq:w} and the one-parameter
group of dilations $\{\varphi_{s}\}_{s\in\rr}$ in
\eqref{eq:dilations}.  Let $\Om_{s}=\varphi_{s}(\Om)$ and
$\Sg_{s}=\ptl\Om_{s}$.  Denote $V(s):=V(\Om_{s})$ and
$A(s):=A(\Sg_{s})$.  From the Christoffel symbols
\eqref{eq:christoffel}, it can be easily proved that $\divv W=4$,
where $\divv W$ is the Riemannian divergence of $W$.  By the first
variation formula of volume \eqref{eq:1stvol} we have
\[
V'(0)=\int_{\Om}\divv W=
4\,V(\Om),
\]
and so $V(s)=e^{4s}V(\Om)$.

Let us calculate now the variation of area $A'(0)$.  Recall that for fixed $s$ and $p\in\hh^1$,
we have $(d\varphi_{s})_{p}(X_{p})=e^s X_{\varphi_{s}(p)}$, $(d\varphi_{s})_{p}(Y_{p})=e^s
Y_{\varphi_{s}(p)}$, and $(d\varphi_{s})_{p}(T_{p})=e^{2s} T_{\varphi_{s}(p)}$.  Let $N$ be the
inner unit normal to $\Sg$, and $p\in\Sg$.  From the calculus of $(d\varphi_{s})_{p}$ we see
that $\varphi_{s}$ preserves the horizontal distribution, so that $p$ lies in the regular part
of $\Sg$ if and only if $\varphi_{s}(p)$ lies in the regular part of $\Sg_s$. Assume $p$ is a
regular point of $\Sg$.  Then we can choose $\alpha$, $\beta\in\rr$ so that $\{e_{1},e_{2}\}$,
with $e_{1}=\cos\alpha\,X_{p}+\sin\alpha\,Y_{p}$, and
$e_{2}=\cos\beta\,(-\sin\alpha\,X_{p}+\cos\alpha\,Y_{p})+\sin\beta\,T_{p}$, is an orthonormal
basis of $T_{p}\Sg$.  For the normal $N$ we have $\pm
N_{p}=-\sin\beta\,(-\sin\alpha\,X_{p}+\cos\alpha\,Y_{p})+\cos\beta\,T_{p}$, and so
$|N_{H}|_{p}=|\sin\beta|$.  We have
$(d\varphi_{s})_{p}(e_{1})=e^s\,(\cos\alpha\,X_{\varphi_{s}(p)}
+\sin\alpha\,Y_{\varphi_{s}(p)})$, and $(d\varphi_{s})_{p}(e_{2})=
e^s\cos\beta\,(-\sin\alpha\,X_{\varphi_{s}(p)}+\cos\alpha\,Y_{\varphi_{s}(p)})
+e^{2s}\sin\beta\,T_{\varphi_{s}(p)}$, and so $|\text{Jac}(\varphi_{s})|_{p}=e^{2s}
(\cos^2\beta+e^{2s}\sin^2\beta)^{1/2}$.  Hence the relation $(d\Sg_s)_{\varphi_{s}(p)}
=e^{2s}(\cos^2\beta+e^{2s}\sin^2\beta)^{1/2}(d\Sg)_{p}$ holds between the area elements of
$\Sg_s$ and $\Sg$.  For the unit normal $N'$ of $\Sg_s$ at $\varphi_{s}(p)$ 
we have 
\begin{align*}
\pm N'_{\varphi_{s}(p)}&=e^{-s}(\cos^2\beta+e^{2s}\sin^2\beta)^{-1/2}
\\
&\times
[-e^{2s}\sin\beta\,(-\sin\alpha\,X_{\varphi_{s}(p)}+
\cos\alpha\,Y_{\varphi_{s}(p)})+e^s\cos\beta\,T_{\varphi_{s}(p)}], 
\end{align*}
and so $|N'_{H}|_{\varphi_{s}(p)} =e^s|\sin\beta|\,(\cos^2\beta+e^{2s}\sin^2\beta)^{-1/2}$.  Hence
\[
|N'_{H}|_{\varphi_{s}(p)}\,(d\Sg_s)_{\varphi_{s}(p)}= e^{3s}|N_{H}|\,(d\Sg)_{p}.
\]
Since $p$ is an arbitrary regular point of $\Sg$, integrating the above displayed formula over
$\Sg-\Sg_0$ and using the area formula we have $A(s)=e^{3s}A(\Sg)$, and so
\[
A'(0)=3\,A(\Sg).
\]
Finally, as $\Sg$ is volume-preserving area-stationary, we deduce from
Corollary~\ref{cor:caract} that $A'(0)=2HV'(0)$, and equality~\eqref{eq:minkowski} follows.
\end{proof}

\begin{corollary}
\label{cor:hpositive} Let $\Sg\subset\hh^1$ be a volume-preserving area-stationary $C^2$
surface enclosing a bounded region $\Om$.  Then the constant mean curvature of the regular part
of $\Sg$ with respect to the inner normal is positive. In particular, there are no compact
area-stationary $C^2$ surfaces in $\hh^1$.
\end{corollary}

\begin{remark}
The generalization of \eqref{eq:minkowski} to the $(2n+1)$-dimensional Heisenberg group $\hh^n$
is immediate. By using the first variation formula in \cite[Lemma 3.2]{revolucion} and the
arguments in this section we get that, for a $C^2$ volume-preserving area-stationary
hypersurface $\Sg\subset\hh^n$ enclosing a bounded region $\Om$, we have
\[
\label{eq:minkowskihn}
(2n+1)\,A(\Sg)=4n(n+1)H\,V(\Om).
\]
\end{remark}

We finish this section with a characterization of area-stationary surfaces in terms of
geometric conditions.  For that, we need additional information on the singular set $\Sg_{0}$
of a constant mean curvature surface $\Sg\sub\hh^1$. The set $\Sg_0$ has been recently studied
by J.-H.~Cheng, J.-F.~Hwang, A.~Malchiodi and P.~Yang \cite{chmy}. Their results are local and
also valid when the mean curvature is bounded on the regular set $\Sg-\Sg_0$.  By Theorem
\ref{th:ruled} we can replace ``characteristic curves'' in their statement by ``geodesics of
the same curvature''. We summarize their results in the following theorem.

\begin{theorem}[{\cite[Theorem~B]{chmy}}]
\label{th:chmy}
Let $\Sg\subset\hh^1$ be a $C^2$ oriented immersed surface with constant mean
curvature $H$.  Then the singular set $\Sg_{0}$ consists of isolated
points and $C^1$ curves with non-vanishing tangent vector.  Moreover,
we have
\begin{itemize}
\item[(i)] $($\cite[Theorem~3.10]{chmy}$)$  If $p\in\Sg_{0}$ is isolated
then there is $r>0$ and $\la\in\rr$ with $|\la|=|H|$ such that the set described as
\[
D_{r}(p)=\{\gamma_{p,v}^\la(s);v\in T_p\Sg,\,|v|=1,\,s\in [0,r)\},
\]
is an open neighborhood of $p$ in $\Sg$. \vspace{0,1cm}
\item[(ii)] $($\cite[Proposition~3.5 and Corollary~3.6]{chmy}$)$ If $p$ is
contained in a $C^1$ curve $\Ga\subset\Sg_{0}$ then there is a neighborhood $B$ of $p$ in
$\Sg$ such that $B-\Gamma$ is the union of two disjoint connected open sets $B^+$ and $B^-$
contained in $\Sg-\Sg_0$, and $\nu_H$ extends continuously to $\Ga$ from both sides of
$B-\Gamma$, i.e., the limits
\[
\nu^+_H(q)=\lim_{x\to q,\ x\in B^+}\nu_H(x),
\qquad
\nu^-_H(q)=\lim_{x\to q,\ x\in B^-}\nu_H(x)
\]
exist for any $q\in\Ga\cap B$. These extensions satisfy $\nu_H^+(q)=-\nu_H^-(q)$. Moreover,
there are exactly two geodesics $\ga_{1} ^\la\sub B^+$ and $\ga_{2}^\la\sub B^-$ starting from
$q$ and meeting transversally $\Ga$ at $q$ with initial velocities
\[
(\ga_{1}^{\la})'(0)=-(\ga^\la_{2})'(0).
\]
The curvature $\la$ does not depend on $q$ and satisfies $|\la|=|H|$.
\end{itemize}
\end{theorem}

\begin{remark}
\label{rem:lambdaorientation}
The relation between $\lambda$ and $H$ depends on the value of the
normal $N$ in the singular point $p$.  If $N_{p}=T$ then $\lambda=H$,
while we have $\lambda=-H$ whenever $N_{p}=-T$.  In case $\la=H$ the
geodesics $\ga^\la$ in Theorem~\ref{th:chmy} are characteristic curves
of $\Sg$.
\end{remark}

In Euclidean space it is equivalent for a surface to be area-stationary (resp.
volume-preserving area-stationary) and to have zero (resp. constant) mean curvature. For a
surface $\Sg$ is $\hh^1$ this also holds if the singular set $\Sg_{0}$ consists only of
isolated points. In the general case, we have the following

\begin{theorem}
\label{th:constant} Let $\Sg\subset\hh^1$ be either an oriented area-stationary $C^2$ immersed
surface or a volume-preserving area-stationary $C^2$ compact surface enclosing a region $\Om$.
Then the mean curvature of\/ $\Sg-\Sg_{0}$ is, respectively, zero or constant and, in both
cases, the characteristic curves meet the singular curves, if they exist, orthogonally. The
converse is also true.
\end{theorem}

\begin{proof}
Suppose first that $\Sg$ is area-stationary.  That the mean curvature is zero or constant on
$\Sg-\Sg_0$ follows from Corollary~\ref{cor:hconstant}. Assume $\Ga$ is a singular curve and
let $p\in\Ga$. By Theorem \ref{th:chmy} (ii) the curve $\Ga$ is $C^1$ and we can take a
neighborhood $B$ of $p$ in $\Sg$ such that $B-\Ga$ consists of the union of two open connected
sets $B^+$ and $B^-$ contained in $\Sg-\Sg_0$.  Let $\xi$ be the unit normal to $\Ga$ in $\Sg$
pointing into $B^+$. Let $f:\Ga\to\rr$ be any $C^1$ function supported on $\Ga\cap B$.  Extend
$f$ to a $C^1$ function $u:B\to\rr$ with compact support in $B$ and mean zero. Since $\Sg$ is
area-stationary, by \eqref{eq:1stvar} and the divergence theorem we have
\begin{align*}
0=A'(0)&=-\int_{B}\divv_{\Sg}\big(u\,(\nuh)^\top\big)\,d\Sg \\
&=
-\int_{B^+}\divv_{\Sg}\big(u\,(\nuh)^\top\big)\,d\Sg
-\int_{B^-}\divv_{\Sg}\big(u\,(\nuh)^\top\big)\,d\Sg \\
&=\int_{\Ga} f\,\escpr{\xi,\nuh^+}\,d\Ga-\int_{\Ga}
f\,\escpr{\xi,\nuh^-}\,d\Ga\\
&=2\int_{\Ga} f\,\escpr{\xi,\nuh^+}\,d\Ga,
\end{align*}
since the extensions $\nuh^+$, $\nuh^-$ of $\nuh$ given in
Theorem~\ref{th:chmy} (ii) satisfy $\nuh^+=-\nuh^-$.  As $f$ is an
arbitrary function on $\Ga\cap B$ we conclude that
$\escpr{\xi,\nuh^+}\equiv 0$ on $\Ga\cap B$.  This means that $\nuh^+$
is tangent to $\Ga\cap B$ and so the two characteristic curves approaching
$p$ meet the singular curve $\Ga$ in an orthogonal way.

We will see the converse for constant mean curvature.  Let $U$ be a $C^1$ vector field inducing
a volume-preserving variation of $\Sg$.  Let $u=\escpr{U,N}$.  By the first variation of volume
\eqref{eq:1stvol} we have $\int_{\Sg} u\,d\Sg=0$.  By \eqref{eq:1stvar}
\[
A'(0)= -\int_\Sg\divv_\Sg\big(u\,(\nu_H)^\top\big)\,d\Sg,
\]
since $u$ has mean zero and $\divv_{\Sg}\nuh$ is a constant.  To
analyze the above integral, we consider disjoint open balls
$B_{\eps}(p_{i})$ (for the Riemannian
distance on $\Sg$) of small radius $\eps>0$, centered at the isolated
points $p_{1},\ldots,p_{k}$ of the singular set $\Sg_{0}$.  By the
divergence theorem in $\Sg$, and the fact that the characteristic
curves meet orthogonally the singular curves we have, for
$\Sg_{\eps}=\Sg-\bigcup_{i=1}^k B_{\eps}(p_{i})$,
\begin{equation*}
-\int_{\Sg_{\eps}}\divv_\Sg\big(u\,(\nu_H)^\top\big)\,d\Sg= \sum_{i=1}^k\, \int_{\ptl
B_{\eps}(p_{i})} u\,\escpr{\xi_{i},\nuh}\,dl,
\end{equation*}
where $\xi_{i}$ is the inner unit normal vector to $\ptl
B_{\eps}(p_{i})$ in $\Sg$. Note also that
\[
\bigg|\sum_{i=1}^k\,\int_{\ptl B_{\eps}(p_{i})}u\,\escpr{\xi_{i},\nu_{H}}\,dl\bigg|
\leq\big(\sup_{\Sg}\,|u|\big)\,\sum_{i=1}^kL(\ptl B_{\eps}(p_{i})),
\]
where $L(\ptl B_{\eps}(p_{i}))$ is the Riemannian length of $\ptl B_{\eps}(p_{i})$. Finally, as
$|\divv_{\Sg}(u\,(\nu_{H}^\top))|\leq
(\sup_{\Sg}\,|u|)\,|\divv_{\Sg}\nu_{H}-|N_{H}|\,\divv_{\Sg}N|+|\nabla_{\Sg}u|\in L^1(\Sg)$, we
can apply the dominated convergence theorem and the fact that $L(\ptl B_{\eps}(p_{i}))\to 0$
when $\eps\to 0$ to prove the claim.
\end{proof}

\begin{example}
Any sphere $\sph_{\la}$ is a volume-preserving area-stationary surface by Theorem
\ref{th:constant} since it has constant mean curvature in $\Sg-\Sg_0$ and $\Sg_0$ consists of
isolated points.
\end{example}

\begin{remark}
Recently, J.-H.~Cheng, J.-F.~Hwang and P.~Yang \cite[Theorem~6.3 and (7.2)]{chy} have obtained
Theorem~\ref{th:constant} when $\Sg$ is a $C^2$ graph over a bounded set $D$ of the $xy$-plane
which is a weak solution of the equation $\divv_\Sg\nuh=-2H$ (\cite[Equation (3.12)]{chy}). As
it is proved in \cite[Theorem~3.3]{chy} such a graph minimizes the functional $A-2HV$ amongst
all graphs $\Sg'$ in the Sobolev space $W^{1,1}(D)$ with $\ptl\Sg'=\ptl\Sg$. In particular,
these graphs are area-stationary for variations by graphs leaving invariant $\ptl\Sg$.
\end{remark}

For a $C^2$ area-stationary surface we can use Theorem \ref{th:constant} to
improve the $C^1$ regularity of  the singular curves obtained in
\cite[Theorem~3.3]{chmy}.

\begin{proposition}
\label{prop:c2singcurves}
If $\Sg$ is a $C^2$ oriented immersed area-stationary surface $($with or
without a volume constraint \!$)$ then any singular curve of $\Sg$ is
a $C^2$ smooth curve.
\end{proposition}

\begin{proof}
By Corollary~\ref{cor:hconstant} we know that $\Sg-\Sg_0$ has constant
mean curvature $H$.  Let $\Ga$ be a connected singular curve of $\Sg$
and $p_{0} \in\Ga$.  By taking the opposite unit normal to $\Sg$ if
necessary we can assume that $N=-T$ along $\Ga$.  By using
Theorem~\ref{th:constant} (ii) and the remark below, we can find a
small neighborhood $B$ of $p_{0}$ in $\Sg$ such that $B^+$ is foliated by
geodesics of the same curvature $\la=H$ reaching $\Ga\cap B$ at
finite, positive time.  These geodesics are characteristic curves 
of $\Sg$ and meet $\Ga$ orthogonally by Theorem \ref{th:constant}.

Let $Z$ be the characteristic vector field of $\Sg$ with respect to
$N$.  Take a point $q\in B^+$ such that $\ga^\la_{q,Z(q)}(s(q))=p_{0}$
for some $s(q)>0$.  We consider a $C^2$ curve $\mathcal{C}\sub B^+$
passing through $q$ and meeting transversally the geodesics only at
one point.  We define the $C^1$ map $F:\mathcal{C}\times
(0,+\infty)\to\hh^1$ given by $F(x,s)=\gamma^\la_{x,Z(x)}(s)$.  For
any $x\in\mathcal{C}$ there is a first value $s(x)>0$ such that
$F(x,s(x))\in\Ga$.  Moreover, by using the orthogonality condition in
Theorem~\ref{th:constant} we can choose the curve $\mathcal{C}$ so
that the differential of $F$ has rank two for any $(x,s(x))$ near to
$(q,s(q))$.  Thus, for some $\delta>0$ we have that
$\Sg'=\{F(x,s);x\in [q-\delta,q+\delta], \ s\in [0,s(x)+\delta]\}$ is
a $C^1$ extension of $\Sg$ beyond the singular curve $\Ga$.  In
particular $\Sg$ and $\Sg'$ are tangent along $\Ga$.  The horizontal
tangent vector to $\Sg'$ given by $Z'=(\ptl F/\ptl
s)(x,s)=(\ga^\la_{x,Z(x)})'(s)$ is a $C^1$ extension of $Z$.  Finally
the orthogonality condition implies that the restriction of $J(Z')$ is
a unit $C^1$ tangent vector to $\Ga$.  We conclude that $\Ga$ is a
$C^2$ smooth curve around $p_{0}$ and the proof follows.
\end{proof}

\section{\textbf{Entire area-stationary graphs in $\hh^1$}}
\label{sec:bernstein}
\setcounter{equation}{0}

An \emph{entire graph} over a plane is one defined over the whole plane.  A classical theorem
by Bernstein shows that the only entire minimal graphs in Euclidean space $\rr^3$ are the
planes.  In \cite[Theorem~D]{pauls}, S. Pauls observed the existence of entire graphs with
$H=0$ in $\hh^1$ different from Euclidean planes.  These are obtained by rotations about the
$t$-axis of a graph of the form
\begin{equation}
\label{eq:bernstein}
t=xy+g(y),\qquad\text{where } \ g\in C^2(\rr).
\end{equation}
In \cite[Theorem~A]{chmy}, J.-H.~Cheng, J.-F.~Hwang, A.~Malchiodi and P.~Yang proved that
Euclidean planes and vertical rotations of \eqref{eq:bernstein} are the unique $C^2$ graphs
over the $xy$-plane with $H=0$, see also \cite[Theorem~D]{gp}.  Here we show that according to
Theorem~\ref{th:constant} not all the graphs in \eqref{eq:bernstein} are area-stationary.  In
precise terms, we have

\begin{theorem}
\label{th:bernstein}
The unique entire $C^2$ area-stationary graphs over the $xy$-plane in
$\hh^1$ are Euclidean planes and vertical rotations of graphs of the form
\[
t=xy+(ay+b),
\]
where $a$ and $b$ are real constants.
\end{theorem}

\begin{proof}
Let $\Sg$ be a $C^2$ entire area-stationary graph over the $xy$-plane in $\hh^1$.  By
Theorem~\ref{th:constant} we know that the mean curvature of $\Sg-\Sg_0$ vanishes and the
intersection between characteristic lines and singular curves is orthogonal. By the
classification in \cite[Theorem~A]{chmy} for entire graphs with $H=0$ we have that $\Sg$ is a
Euclidean plane or a vertical rotation of \eqref{eq:bernstein}.  That Euclidean planes are
area-stationary follows from Theorem~\ref{th:constant} since they only have isolated
singularities. To prove the claim we suppose that $\Sg$ coincides with \eqref{eq:bernstein}.
The surface $\Sg$ has a connected curve $\Gamma$ of singular points whose projection to the
$xy$-plane is given by the equation $x=-g'(y)/2$. We can parameterize $\Ga$ by
\[
\Ga(s)=\left(-\frac{g'(s)}{2},s,g(s)-\frac{g'(s)\,s}{2}\right),\qquad s\in\rr,
\]
and so, if $\Ga(s_{0})=p_{0}$, then $\dot{\Ga}(s_{0})=(-g''(s_{0})/2)\,X_{p_{0}}+Y_{p_{0}}$. On
the other hand, it is not difficult to check that for a fixed $y\in\rr$, the straight line
$t=xy+g(y)$ is a characteristic curve of $\Sg$ when removing the contact point with $\Ga$.  We
parameterize this line as
\[
S_{y}(s)=(s,y,sy+g(y)),\qquad s\in\rr,
\]
so that if $S_{y}(s_{1})=p_{0}$ then $\dot{S}_{y}(s_{1})=X_{p_{0}}$.  From these computations
we see that, for $p_{0}=\Ga(s_{0})=S_{y}(s_{1})$ we have
\[
\escpr{\dot{\Ga}(s_0),\dot{S_y}(s_1)}=-\frac{g''(y)}{2}.
\]
We conclude that the characteristic lines $S_y$ meet orthogonally the singular curve $\Ga$ if
and only if $g(y)=ay+b$ for some real constants $a$ and $b$.
\end{proof}

\begin{remark}
While Euclidean planes have only an isolated singular point, the
entire area-stationary graphs obtained by rotations of $t=xy+(ay+b)$
have a straight line of singular points.  From a geometric point of
view, these second surfaces are constructed by taking a horizontal
straight line $R$ and attaching at any point of $R$ the unique
straight line which is both horizontal and orthogonal to $R$.  The
remaining surfaces defined by equation \eqref{eq:bernstein} have
vanishing mean curvature outside the singular set, but they are not
area-stationary.
\end{remark}

We finish this section showing that the graphs obtained in
Theorem~\ref{th:bernstein} are globally area-minimizing.  This is a
counterpart in $\hh^1$ of a well-known result for minimal graphs in
$\rr^3$.

We say that a surface $\Sg\sub\hh^1$ is \emph{area-minimizing} if any region $M\sub\Sg$ has
less area than any other $C^1$ compact surface $M'$ in $\hh^1$ with $\ptl M=\ptl M'$.  In
\cite[Proposition~6.2]{chmy} it was proved by using a calibration argument that any $C^2$
surface in $\hh^1$ with vanishing mean curvature locally minimizes the area around any point in
the regular set.  Here, we adapt the calibration argument in order to deal with surfaces with
singularities, and we obtain

\begin{theorem}
\label{th:areaminimizing}
Any entire $C^2$ area-stationary graph $\Sg$ over the $xy$-plane in
$\hh^1$ is area-minimizing.
\end{theorem}

\begin{proof}
After a vertical rotation about the $t$-axis we may assume, by
Theorem~\ref{th:bernstein}, that $\Sg$ coincides with a Euclidean
plane or with a graph of the form $t=xy+ay+b$, for some $a,b\in\rr$.
Let $\Sg_{t}$ be area-stationary graph obtained by applying to $\Sg$
the left translation $L_{t}$ by the vertical vector $tT$.  The family
$\{\Sg_{t}\}_{t\in\rr}$ is a foliation of $\hh^1$ by area-stationary
surfaces.  Moreover, $L_{t}$ preserves the horizontal distribution and
hence $p\in\Sg-\Sg_{0}$ if and only if
$L_{t}(p)\in\Sg_{t}-(\Sg_{t})_{0}$.  Therefore, the set $P=\bigcup_{t}
(\Sg_{t})_{0}$ is either a vertical straight line if $\Sg$ is a plane
or a vertical plane if $\Sg$ is a graph $t=xy+ay+b$.  Consider a $C^1$
vector field $N$ on $\hh^1$ so that the restriction $N_{t}$ of $N$ to 
$\Sg_{t}$ is a unit normal vector to $\Sg_{t}$. We denote
$N_{H}/|N_{H}|$ by $\nuh$, and $Z=J(\nuh)$, which are $C^1$ vector
fields on $\hh^1-P$. 

Let us compute $\divv\nuh$.  Take a point $p$ in the regular set of
$\Sg_{t}$ for some $t\in\rr$.  We have an orthonormal basis of
$T_{p}\hh^1$ given by $\{Z_{p},(\nuh)_{p},T\}$.  Denote by $H_{t}$ the
mean curvature of $\Sg_{t}$ with respect to $N_{t}$.  By using
equation \eqref{eq:dznuh} and that $\nuh$ is a horizontal unit vector
field, we get
\begin{align*}
\divv\nuh&=\escpr{D_{Z}\nuh,Z}+\escpr{D_{\nuh}\nuh,\nuh}+
\escpr{D_{T}\nuh,T}
\\
&=-2H_{t}-\escpr{\nuh,D_{T}T}=0,
\end{align*}
where in the last equality we have used that $H_{t}\equiv 0$ since
$\Sg_{t}$ is area-stationary (Corollary~\ref{cor:hconstant}(i)), and that
$D_{T}T=0$.

Consider a region $M\sub\Sg$ and a compact $C^1$ surface $M'\sub\hh^1$
with $\ptl M=\ptl M'$.  We denote by $\Om$ the open set bounded by $M$
and $M'$.  The set $\Om$ has finite perimeter in the Riemannian
manifold $(\hh^1,g)$ since it is bounded and the two-dimensional
Riemannian Hausdorff measure of $\ptl\Om\cap C$ is finite for any
compact set $C\sub\hh^1$, see \cite[Theorem~1, p.~ 222]{evans}.  For the
following arguments we may assume $\Om$ connected, and that
$\ptl\Om=M\cup M'$.  We fix the outward normal vector $N$ to $\Sg$,
and the unit normal vector $N'$ to $M'$, to point into $\Om$.  As a
consequence, we can apply the Gauss-Green Theorem for sets of finite
perimeter \cite[Theorem~1, p.~209]{evans} so that, for any $C^1$ vector
field $U$ on $\hh^1$, we have
\begin{equation}
\label{eq:gaussgreen}
\int_{\Om}\divv U\,dv=\int_{M}\escpr{U,N}\,dM-\int_{M'}\escpr{U,N'}\,dM'.
\end{equation}

In order to prove $A(M)\leq A(M')$ we distinguish two cases.

Case 1.  If $\Sg$ is a Euclidean plane, then $\nuh$ is defined in the
closure of $\Om$ outside a set contained in a straight line.  Thus, we
can apply \eqref{eq:gaussgreen} to deduce
\begin{align*}
0=\int_{\Om}\divv\nuh\,dv&=\int_{M}\escpr{\nuh,N}\,dM
-\int_{M'}\escpr{\nuh,N'}\,dM'
\\
&=\int_{M}|N_{H}|\,dM-\int_{M'}\escpr{\nuh,N'_{H}}\,dM'
\\
&\geq A(M)-A(M').
\end{align*}
To obtain the last inequality we have used the Cauchy-Schwarz
inequality and that $|\nuh|=1$.  This proves the claim.

Case 2.  If $\Sg$ is a graph of the form $t=xy+ay+b$, then $\nuh$ is
defined on $\Om-P$, where $P$ is a vertical Euclidean plane.  Denote
by $P^+$ and $P^-$ the open half-planes determined by $P$.  For any
set $E\sub\hh^1$, we let $E^+=E\cap P^+$ and $E^-=E\cap P^-$.  The
sets $\Om^+$ and $\Om^-$ has finite perimeter in $(\hh^1,g)$.
Moreover, by Theorem \ref{th:chmy} (ii) the vector field $\nuh$ extends 
continuously to $P$ from $\Om^+$ and $\Om^-$ .  Therefore
\begin{align*}
0&=\int_{\Om^+}\divv\nuh\,dv=\int_{M^+}\escpr{\nuh,N}\,dM
-\int_{(M')^+}\escpr{\nuh,N'}\,dM'-\int_{\Om\cap P}\escpr{\nuh^+,\xi}\,dP
\\
0&=\int_{\Om^-}\divv\nuh\,dv=\int_{M^-}\escpr{\nuh,N}\,dM
-\int_{(M')^-}\escpr{\nuh,N'}\,dM'+\int_{\Om\cap P}\escpr{\nuh^-,\xi}\,dP,
\end{align*}
where $\xi$ is the unit normal vector to $P$ pointing into $\Om^+$.
As $\nuh^+=-\nuh^-$, by summing the previous equalities we deduce
\begin{align*}
0&=\int_{M}\escpr{\nuh,N}\,dM-\int_{M'}\escpr{\nuh,N'}\,dM' -2\int_{\Om\cap
P}\escpr{\nuh^+,\xi}\,dP
\\
&\geq A(M)-A(M')-2\int_{\Om\cap P}\escpr{\nuh^+,\xi}\,dP.
\end{align*}
Finally, the orthogonality condition between characteristic lines and
singular curves in Theorem \ref{th:constant} implies that 
$\escpr{\nuh^+,\xi}=0$ on $\Om\cap P$.  Thus, we get $A(M)\leq~A(M')$.
\end{proof}

\begin{remark}
If $\Sg$ is an area-stationary surface in $\hh^1$, and there is a left
invariant vector field $V$ in $\hh^1$ transverse to $\Sg$, then we can
produce a local foliation by area-stationary surfaces around $\Sg$ by
using the flow associated to $V$.  The arguments in the proof of
Theorem~\ref{th:areaminimizing} show that $\Sg$ is locally
area-minimizing, i.e., bounded portions of $\Sg$ minimize area amongst
surfaces with boundary on $\Sg$ and contained in the foliated
neighborhood of $\Sg$.
\end{remark}

\begin{remark}
1. It follows from \cite[Proposition~6.2 and Theorem~3.3]{chy} that a $C^2$ 
area-stationary graph over a bounded domain $D$ of the $xy$-plane 
minimizes the area amongst all graphs $\Sg'$ in the Sobolev space 
$W^{1,1}(D)$ with $\ptl\Sg'=\ptl\Sg$. This has been recently improved in 
\cite[Example 2.7]{bscv} where it is shown that such a graph is 
area-minimizing.    

2.  Theorem~\ref{th:areaminimizing} does not hold for a graph over the
$xt$-plane, see an example in \cite{dgn3}.  In
\cite[Theorem~5.3]{bscv} it is proved that the unique $C^2$ entire,
area-minimizing \emph{intrinsic graphs} over the $xt$-plane are
vertical planes.
\end{remark}

\section{\textbf{Complete volume-preserving area-stationary
surfaces in $\hh^1$}}
\label{sec:mainresult}
\setcounter{equation}{0}

An immersed surface $\Sg\subset\hh^1$ is \emph{complete} if it is
complete in the Riemannian manifold $(\hh^1,g)$.  Completeness for a
constant mean curvature surface is equivalent to that the singular
curves in $\Sg_{0}$ are closed in $\hh^1$ and that characteristic
curves in $\Sg-\Sg_{0}$ extend up to singular points of $\Sg$.

In this section we obtain classification results for complete
area-stationary surfaces under a volume constraint in $\hh^1$.  We say
that a complete noncompact oriented $C^2$ surface in $\hh^1$ is
volume-preserving area-stationary if it has constant mean curvature
off of the singular set and the characteristic curves meet
orthogonally the singular curves.  By Theorem \ref{th:constant} this
implies that the surface is area-stationary for any variation with
compact support of the surface such that the volume
\eqref{eq:immersedvolume} of the perturbed region remains constant.

We begin with the description of constant mean curvature surfaces with
isolated singularities.  It was shown in \cite[Proof of Theorem~A]{chmy}
(see also \cite[Proposition~2.1]{ch}) and \cite[Lemma~8.2]{gp} that any
$C^2$ surface with vanishing mean curvature and an isolated singular
point must coincide with a Euclidean plane.  By using the local
behavior of a constant mean curvature surface around a singular point
(Theorem~\ref{th:chmy}) we can prove the following

\begin{theorem}
\label{th:spheres}
Let $\Sg$ be a complete, connected, $C^2$ oriented immersed surface in
$\hh^1$ with non-vanishing constant mean curvature.  If $\Sg$ contains
an isolated singular point then $\Sg$ is congruent with a sphere
$\mathbb{S}_{H}$.
\end{theorem}

\begin{proof}
We choose the unit normal $N$ to $\Sg$ such that the mean curvature
$H$ is positive.  Let $p$ be an isolated singular point of $\Sg$.  By
applying to $\Sg$ the left translation $(L_{p})^{-1}$ we can assume
that $p=0$ and the tangent plane $T_{p}\Sg$ coincides with the
$xy$-plane.  Suppose that $N_p=T$.  For any $r>0$ we consider the set
\[
D_{r}=\{\ga^H_{0,v}(s);\,|v|=1, \ s\in [0,r)\}.
\]
It is clear that the union of $D_{r}$, for $r\in (0,\pi/H)$, coincides
with the sphere $\sph_{H}$ removing the north pole (see
Example~\ref{ex:spheres}).  By Theorem~\ref{th:chmy} (i) and
Remark~\ref{rem:lambdaorientation}, we can find $r>0$ such that
$D_{r}\sub\Sg$.  Let $R=\sup\,\{r>0\ ;\ D_{r}\sub\Sg\}$.  As $\Sg$ is
complete and connected, and $\sph_{H}$ is compact, to prove the claim
it suffices to see that $R=\pi/H$.

Suppose that $R<\pi/H$.  In this case we would have
$\overline{D}_R\sub\Sg$ and so, $\Sg$ and $\sph_{H}$ would be tangent
along the curve $\ptl D_{R}$.  In particular, this curve is contained
in the regular set of $\Sg$.  By Theorem \ref{th:ruled} the
characteristic curve of $\Sg$ passing through any $q\in\ptl D_{R}$ is
an open arc of a geodesic of curvature $H$.  By the uniqueness of the
geodesics this would imply that we may extend any $\ga^H_{0,v}$ inside
$\Sg$ beyond $\ptl D_{R}$, a contradiction with the definition of $R$.
This proves $R\geq\pi/H$.  On the other hand, $R>\pi/H$ would imply
that $\Sg$ contains a neighborhood of a tangent point between two
different spheres of the same curvature which is not possible since
$\Sg$ is immersed.

Finally, if $N_p=-T$ we repeat the previous arguments by using
geodesics of curvature $-H$ and we conclude that $\Sg$ coincides with
a vertical translation of $\sph_H$.
\end{proof}

Theorem~\ref{th:spheres} does not provide information about
non-vanishing constant mean curvature surfaces in $\hh^1$ with at
least one singular curve.  We will treat this situation in the
particular case of volume-preserving area-stationary surfaces, where
we have by Theorem~\ref{th:constant} the additional condition that the
characteristic curves meet orthogonally the singular curves.  We first
study in more detail the behavior of the characteristic curves far
away from a singular curve.

Let $\Ga$ be a $C^2$ horizontal curve in $\hh^1$.  We parameterize
$\Ga=(x,y,t)$ by arc-length $\eps\in I$, where $I$ is an open
interval.  The projection $\alpha=(x,y)$ is a plane curve with
$|\dot{\alpha}|=1$.  We denote by $h$ the planar geodesic curvature
of $\alpha$ with respect to the unit normal vector
$(-\dot{y},\dot{x})$, that is $h=\dot{x}\,\ddot{y}-\ddot{x}\,\dot{y}$.
As $\ga$ is horizontal, we have $\dot{t}=\dot{x}y-x\dot{y}$.
Fix $\la\neq 0$.  For any $\eps\in I$ let $\ga_{\eps}$ be the unique
geodesic of curvature $\la$ with initial conditions
$\ga_{\eps}(0)=\Ga(\eps)$ and
$\dot{\ga_{\eps}}(0)=J(\dot{\Ga}(\eps))$.  We consider the family of
all these geodesics orthogonal to $\Ga$ parameterized by
$F(\eps,s)=\ga_{\eps}(s)=(x(\eps,s),y(\eps,s),t(\eps,s))$, for
$\eps\in I$ and $s\in[0,\pi/|\la|]$.  By equation \eqref{eq:geocoor}
we have
\begin{align}
\nonumber
x(\eps,s)&=x(\eps)-\dot{y}(\eps)\,
\bigg(\frac{\sin(2\lambda\,s)}{2\lambda}\bigg)
+\dot{x}(\eps)\,
\bigg(\frac{1-\cos(2\lambda\,s)}{2\lambda}\bigg),
\\
\label{eq:geocoor2}
y(\eps,s)&=y(\eps)+\dot{y}(\eps)\,
\bigg(\frac{1-\cos(2\lambda\,s)}{2\lambda}\bigg)+\dot{x}(\eps)\,
\bigg(\frac{\sin(2\lambda\,s)}{2\lambda}\bigg),
\\
\nonumber
t(\eps,s)&=t(\eps)+\frac{1}{2\lambda}\,
\bigg(s-\frac{\sin(2\lambda\,s)}{2\lambda}\bigg)-
(x(\eps)\,\dot{x}(\eps)+y(\eps)\,\dot{y}(\eps))\,\left(\frac{\sin(2\lambda
s)}{2\la}\right)
\\
\nonumber &\quad +(\dot{x}(\eps)\,y(\eps)-x(\eps)\,\dot{y}(\eps))\,\left(\frac{1-\cos(2\lambda
s)}{2\la}\right).
\end{align}
From the equations above we see that $F$ is a $C^1$ map.  Clearly
$(\ptl F/\ptl s)(\eps,s)=\dot{\ga}_{\eps}(s)$.  We denote
$V_{\eps}(s):=(\ptl F/\ptl\eps)(\eps,s)$.  In the next result we show
some properties of~$V_{\eps}$.

\begin{lemma}
\label{lem:jacobi3}
In the situation above, $V_{\eps}$ is a Jacobi vector field along
$\ga_{\eps}$ with $V_{\eps}(0)=\dot{\Ga}(\eps)$.  For any $\eps\in I$
there is a unique $s_{\eps}\in (0,\pi/|\la|)$ such that
$\escpr{V_{\eps}(s_{\eps}),T}=0$.  We have $\escpr{V_{\eps},T}<0$ on
$(0,s_{\eps})$ and $\escpr{V_{\eps},T}>0$ on $(s_{\eps} ,\pi/|\la|)$.
Moreover $V_{\eps}(s_{\eps})= J(\dot{\ga}_{\eps}(s_{\eps}))$.
\end{lemma}

\begin{proof}
By the definition of $V_{\eps}$ we have $V_{\eps}(0)=\dot{\Ga}(\eps)$
and
\[
V_{\eps}(s)=
\frac{\ptl x}{\ptl\eps}(\eps,s)\,\,X
+\frac{\ptl y}{\ptl\eps}(\eps,s)\,\,Y+ \left(\frac{\ptl
t}{\ptl\eps}-y\,\,\frac{\ptl x}{\ptl\eps}+x\,\,\frac{\ptl
y}{\ptl\eps}\right)(\eps,s)\,\,T.
\]
The Euclidean components of $V_{\eps}(s)$ are easily computed from
\eqref{eq:geocoor2}, so that we obtain
\begin{align*}
\nonumber
\frac{\ptl x}{\ptl\eps}(\eps,s)&=\dot{x}(\eps)-\ddot{y}(\eps)\,
\bigg(\frac{\sin(2\lambda\,s)}{2\lambda}\bigg)
+\ddot{x}(\eps)\,
\bigg(\frac{1-\cos(2\lambda\,s)}{2\lambda}\bigg),
\\
\frac{\ptl y}{\ptl\eps}(\eps,s)&=\dot{y}(\eps)+\ddot{y}(\eps)\,
\bigg(\frac{1-\cos(2\lambda\,s)}{2\lambda}\bigg)+\ddot{x}(\eps)\,
\bigg(\frac{\sin(2\lambda\,s)}{2\lambda}\bigg),
\\
\nonumber
\frac{\ptl t}{\ptl\eps}(\eps,s)&=\dot{t}(\eps)+\frac{1}{2\lambda}\,
\bigg(s-\frac{\sin(2\lambda\,s)}{2\lambda}\bigg)-
(1+x(\eps)\,\ddot{x}(\eps)+y(\eps)\,\ddot{y}(\eps))\,\left(\frac{\sin(2\lambda
s)}{2\la}\right)
\\
\nonumber &\quad +(\ddot{x}(\eps)\,y(\eps)-x(\eps)\,\ddot{y}(\eps))\,\left(\frac{1-\cos(2\lambda
s)}{2\la}\right).
\end{align*}
We deduce that $V_{\eps}$ is
$C^\infty$ vector field along $\ga_{\eps}$ and
\[
\escpr{V_{\eps}(s),T}=\frac{1}{\la}\,\left(\frac{1-\cos(2\la\, s)}{2\la}\,\,
h(\eps)-\sin(2\la s)\right),\qquad s\in [0,\pi/|\la|].
\]
That $V_{\eps}$ is a Jacobi vector field along $\ga_{\eps}$ follows
from Lemma~\ref{lem:jacobi2} since $V_{\eps}$ is associated to a
variation of $\ga_{\eps}$ by geodesics of the same curvature.  On the
other hand, the equation above implies that
$\escpr{V_{\eps}(s_{\eps}),T}=0$ for some $s_{\eps}\in (0,\pi/|\la|)$
if and only if
\begin{equation}
\label{eq:despeje}
h(\eps)=\frac{2\la\,\sin(2\la\,s_{\eps})}{1-\cos(2\la\,s_{\eps})}.
\end{equation}
The existence and uniqueness of $s_{\eps}$, and the sign of
$\escpr{V_{\eps},T}$ are consequences of the fact that the function
$f(x)=\sin(x)\,(1-\cos(x))^{-1}$ is periodic, decreasing on $(0,2\pi)$
and satisfies $\lim_{x\to0^+}f(x)=+\infty$ and $\lim_{x\to
(2\pi)^-}f(x)=-\infty$.

Now we use Lemma \ref{lem:jacobi1} and the fact that
$V_{\eps}(0)=\dot{\Ga}(\eps)$ to deduce that the function
$\la\,\escpr{V_{\eps},T}+\escpr{V_{\eps},\dot{\ga}_{\eps}}$ vanishes
along $\ga_{\eps}$. In particular, $V_{\eps}(s_{\eps})$ is a
horizontal vector orthogonal to $\dot{\ga}_{\eps}(s_{\eps})$. Finally,
we have, for $s\in[0,\pi/|\la|]$,
\[
\escpr{V_{\eps}(s),J(\dot{\ga}_{\eps}(s))}=
\bigg(
-\frac{\ptl x}{\ptl\eps}\,\frac{\ptl y}{\ptl s}
+\frac{\ptl y}{\ptl\eps}\,\frac{\ptl x}{\ptl s}
\bigg)(\eps,s)=
\frac{\sin(2\la\,s)}{2\la}
\,h(\eps)-\cos(2\la\,s),
\]
which is equal to $1$ for $s=s_{\eps}$ by \eqref{eq:despeje}.
\end{proof}

The following proposition provides a method to construct immersed
surfaces with constant mean curvature in $\hh^1$ bounded by two
singular curves.  Geometrically we only have to leave from a given
horizontal curve by segments of orthogonal geodesics of the same
curvature.  The length of these segments depends on the \emph{cut
function} $s_{\eps}$ introduced in Lemma~\ref{lem:jacobi3}.  We also
characterize when the resulting surface is volume-preserving
area-stationary.

\begin{proposition}
\label{prop:sigmala}
Let $\Ga$ be a $C^{k+1}$ $(k\geq 1)$ horizontal curve in $\hh^1$
parameterized by arc-length $\eps\in I$.  Consider the map
$F(\eps,s)=\ga_{\eps}(s)$, where $\ga_{\eps}:[0,\pi/|\la|]\to\hh^1$ is
the geodesic of curvature $\la\neq 0$ with initial conditions
$\Ga(\eps)$ and $J(\dot{\Ga}(\eps))$.  Let $s_{\eps}$ be the function 
introduced in Lemma~\ref{lem:jacobi3}, and let
$\Sg_{\la}(\Ga)=\{F(\eps,s);\,\eps\in I, \,s\in [0,s_{\eps}]\}$.  Then
we have
\begin{itemize}
 \item[(i)] $\Sg_{\la}(\Ga)$ is an immersed surface of class $C^k$ in $\hh^1$.
 \item[(ii)] The singular set of $\Sg_\la(\Ga)$ consists of two curves
 $\Ga(\eps)$ and $\Ga_{1} (\eps)=F(\eps,s_{\eps})$.
 \item[(iii)] There is a $C^{k-1}$ unit normal vector $N$ to $\Sg_{\la}(\Ga)$
 such that $N=T$ on $\Ga$ and $N=-T$ on $\Ga_{1}$.
 \item[(iv)] Any $\ga_{\eps}:(0,s_{\eps})\to\hh^1$ is a characteristic curve of
 $\Sg_{\la}(\Ga)$.  In particular, if $k\geq 2$ then $\Sg_{\la}(\Ga)$
 has constant mean curvature $\la$ with respect to $N$.
 \item[(v)] If\/ $\Ga_{1}$ is a $C^2$ smooth curve then the geodesics
 $\ga_{\eps}$ meet orthogonally $\Ga_{1}$ if and only if 
 $s_{\eps}$ is constant along $\Ga$.  This is equivalent to that the
 $xy$-projection of $\Ga$ is either a line segment or a piece of a
 planar circle.
\end{itemize}
\end{proposition}

\begin{proof}
As $\Ga$ is $C^{k+1}$ and the geodesics $\ga_\eps$ depend $C^1$
smoothly on the initial conditions we get that $F$ is a map of class
$C^k$.  Let us consider the vector fields $(\ptl
F/\ptl\eps)(\eps,s)=V_{\eps}(s)$ and $(\ptl F/\ptl
s)(\eps,s)=\dot{\ga}_{\eps}(s)$.  By using Lemma \ref{lem:jacobi3} we
deduce that the differential of $F$ has rank two for any $(s,\eps)\in
I\times [0,\pi/|\la|)$, and that the tangent plane to $\Sg_{\la}(\Ga)$
is horizontal only for the points in $\Ga$ and $\Ga_{1}$.  This proves
(i) and (ii).

Now define the $C^{k-1}$ unit normal vector to the immersion
$F:I\times [0,\pi/|\la|)\to\hh^1$ given by
$N(\eps,s)=|V_{\eps}(s)\wedge\dot{\ga_{\eps}}(s)|^{-1}\,
(V_{\eps}(s)\wedge\dot{\ga_{\eps}}(s))$.  To compute $N$ along $\Ga$
and $\Ga_{1}$ it suffices to use $v\wedge J(v)=T$ for any unit
horizontal vector $v$ together with the fact that
$V_{\eps}(0)=\dot{\Ga}(\eps)$ and
$V_{\eps}(s_{\eps})=J(\dot{\ga}_{\eps}(s_{\eps}))$.  It is easy to see
that the characteristic vector field $Z$ to the immersion is given by
\[
Z(\eps,s)=-\frac{\escpr{V_{\eps}(s),T}}{|\escpr{V_{\eps}(s),T}|}\,\,
\dot{\ga_{\eps}}(s),\qquad\eps\in I, \ \ s\neq 0, s_{\eps}.
\]
By using Lemma \ref{lem:jacobi3} it follows that $Z(\eps,s)=
\dot{\ga_{\eps}}(s)$ whenever $s\in (0,s_{\eps})$. This fact and
Theorem~\ref{th:ruled} prove (iv).

Finally, suppose that $\Ga_{1}$ is a $C^2$ smooth curve (which is
immediate is $k\ge 3$).  The cut function $s(\eps)=s_{\eps}$ is $C^1$
since the graph $(\eps,s(\eps))$ coincides, up to the $C^1$ immersion
$F$, with $\Ga_1$.  The tangent vector to $\Ga_{1}$ is given by
\[
\dot{\Ga}_{1}(\eps)=V_{\eps}(s_{\eps})+\dot{s}(\eps)\,
\dot{\ga_{\eps}}(s_{\eps}).
\]
As $V_\eps(s_{\eps})=J(\dot{\ga}_{\eps}(s_{\eps}))$, we conclude that
the geodesics $\ga_\eps$ meet $\Ga_1$ orthogonally if and only if
$s(\eps)$ is a constant function.  As a consequence, we deduce from
\eqref{eq:despeje} that the planar geodesic curvature of the
$xy$-projection of $\Ga$ is constant and so, this plane curve must
coincide with a line segment or a piece of a Euclidean circle.
\end{proof}

\begin{remark}
\label{re:reverse} 1.  In the proof above it is shown that if we
extend $\Sg_\la(\Ga)$ by the geodesics $\ga_\eps$ beyond the singular
curve $\Ga_1$ then the resulting surface has mean curvature $-\la$
beyond $\Ga_{1}$.  As indicated in Theorem~\ref{th:chmy} (ii), in
order to get an extension of $\Sg_{\la}(\Ga)$ with constant mean
curvature $\la$ we must leave from $\Ga_{1}$ by geodesics of curvature
$-\lambda$ (we must arrive at $\Ga_{1}$ by geodesics of curvature $\la$).

2.  The singular curves $\Ga$ and $\Ga_{1}$ of the surface
$\Sg_{\la}(\Ga)$ could coincide.  We will illustrate this situation
in Example \ref{ex:helices}.
\end{remark}

\begin{remark}
\label{re:reverse2}
Let $\Ga$ be a $C^{k+1}$ $(k\geq 1)$ horizontal curve in $\hh^1$
parameterized by arc-length $\eps\in I$.  We consider the family of
geodesics $\widetilde{\ga}_{\eps}:[0,\pi/|\la|]\to\hh^1$ with
curvature $\la\neq 0$ and initial conditions $\Ga(\eps)$ and
$-J(\dot{\Ga}(\eps))$.  By following the arguments in Lemma
\ref{lem:jacobi3} and Proposition \ref{prop:sigmala} we can construct
the surface
\[
\widetilde{\Sg}_{\la}(\Ga):=\{\widetilde{\ga}_{\eps}(s)\ ;\
\eps\in I,\, s\in [0,s_{\eps}]\},
\]
which is bounded by two singular curves $\Ga$ and $\Ga_2$.  The cut
function $\widetilde{s}_{\eps}$ associated to $\Ga_{2}$ is defined by
the equality $\escpr{\widetilde{V}_{\eps}(\widetilde{s}_{\eps}),
T}=0$, where $\widetilde{V}_{\eps}$ is the Jacobi vector field
associated to $\{\widetilde{\ga}_{\eps}\}$. It is easy to see that
$\widetilde{s}_{\eps}$ satisfies
\[
h(\eps)=\frac{-2\la\sin(2\la\widetilde{s}_{\eps})}{1-\cos(2\la\widetilde{s}_{\eps})}.
\]
From \eqref{eq:despeje} it follows that $s_{\eps}+\widetilde{s}_{\eps}
=\pi/|\la|$.  The vector field $\widetilde{V}_{\eps}$ coincides
with $-J(\dot{\widetilde{\ga}}_{\eps})$ for $s=\widetilde{s}_{\eps}$.
The unit normal $\widetilde{N}$ to $\widetilde{\Sg}_{\la}(\Ga)$ equals
$T$ on $\Ga$ and $-T$ on $\Ga_2$.  When $k\ge 2$, we deduce that the union of
$\Sg_{\la}(\Ga)$ and $\widetilde{\Sg}_{\la}(\Ga)$ is an oriented
immersed surface with constant mean curvature $\la$ and at most three
singular curves.
\end{remark}

Now we shall use Proposition \ref{prop:sigmala} and
Remark~\ref{re:reverse2} to obtain new examples of complete
volume-preserving area-stationary surfaces in $\hh^1$ with singular
curves.  We know by Proposition~\ref{prop:sigmala} (iv) that the
$xy$-projection of the initial curve $\Ga$ must be either a straight
line or a planar circle.  We shall consider the two cases.

\begin{figure}[h]
\centering{\includegraphics[height=5cm]{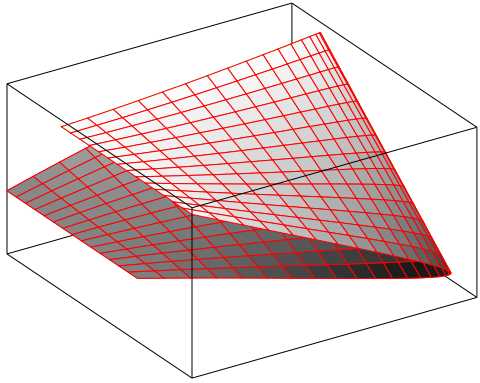}}
\caption{A portion of the surface $\mathcal{S}_{\la}$ composed of
geodesics of curvature $\la>0$ joining two horizontal and parallel
straight lines.}
\label{fig:cylinder1}
\end{figure}

\begin{figure}[h]
\centering{\includegraphics[height=4cm]{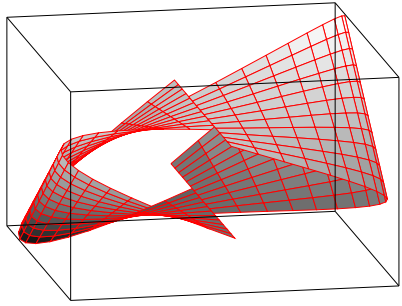}}
\caption{The complete surface $\mathcal{S}_{\la}$.}
\label{fig:cylinder2}
\end{figure}

\begin{example}[Cylindrical surfaces $\mathcal{S}_\la$]
\label{ex:cilindros}
Consider the $x$-axis in $\rr^3$ parameterized by
$\Ga(\eps)=(\eps,0,0)$.  For any $\la\neq 0$ we denote by
$\mathcal{S}_{\la}$ the union of the surfaces $\Sg_{\la}(\Ga)$ and
$\widetilde{\Sg}_{\la}(\Ga)$ constructed in Proposition
\ref{prop:sigmala} and Remark \ref{re:reverse2}.  The surface
$\mathcal{S}_{\la} $ is $C^\infty$ outside the singular curves and has
constant mean curvature $\la$.  The cut functions $s_\eps$ and
$\widetilde{s}_\eps$ can be computed from \eqref{eq:despeje} and the
relation $s_\eps+\widetilde{s}_\eps=\pi/|\la|$, so that, by using
$h_{\eps}\equiv 0$, we get $s_\eps=\widetilde{s}_\eps=\pi/|2\la|$.
From \eqref{eq:geocoor2} we see that the singular curves $\Ga_{1}$ and
$\Ga_{2}$ are different parameterizations of the same curve, namely,
the $x$-axis translated
by the vertical vector $(\text{sgn}(\la)\,\pi/(4\la^2))\,T$, where
$\text{sgn}(x)$ is the sign of $x\in\rr$.  A straightforward
computation from \eqref{eq:geocoor2} shows that $\mathcal{S}_{\la}$ is
the union of the graphs of the functions $f$ and $g$ defined on the
$xy$-strip $-1/|2\la|\leq y\leq 1/|2\la|$ by
\begin{align*}
f(x,y)&=\frac{\text{sgn}(y)}{2\la}\,\left (\displaystyle\frac{\arcsin
    (2\la\,y)}{2\la}-y\,\sqrt{1-4\la\,y^2}\right)-xy,
    \\
g(x,y)&=\frac{1}{2\la}\,\left(\frac{\text{sgn}(\la)\,\pi-\text{sgn}(y)\,
\arcsin(2\la\,y)}{2\la}+\text{sgn}(y)\,y\,\sqrt{1-4\la^2\,y^2}\right)-xy.
\end{align*}
Both functions coincide on the boundary of the strip.  Moreover, it is
easy to see that $\mathcal{S}_{\la}$ is $C^2$ smooth around $\Ga$ and
$\Ga_1=\Ga_2$ but not $C^3$ since
\[
\frac{\ptl^3 f}{\ptl y^3}\,(x,y)=-\frac{\ptl^3 g}{\ptl y^3}\,(x,y)=\text{sgn}(y)\,\,
\frac{8\la\,(1+2\la^2\,y^2)}{(1-4\la^2\,y^2)^{5/2}}.
\]
Finally, an easy argument proves that
$\text{sgn}(\la)\,f(x,y)<\text{sgn}(\la)\,g(x,y)$
for any $(x,y)$ such that $-1/|2\la|<y<1/|2\la|$.  We conclude
that $\mathcal{S}_\la$ is a complete volume-preserving area-stationary
embedded cylinder in $\hh^1$ with two singular curves given by
parallel straight lines, see Figures~\ref{fig:cylinder1} and
\ref{fig:cylinder2}.
\end{example}

\begin{example}[Helicoidal surfaces $\mathcal{L}_\la$]
\label{ex:helices}
Let $\Ga$ be the helix of radius $r>0$ and pitch $\pi/(2r^2)$ in $\rr^3$
given by
\[
\Ga(\eps)=\left(\frac{\sin(2r\eps)}{2
r},\frac{\cos(2r\eps)-1}{2r},\frac{1}{2r}\,
\left(\eps-\frac{\sin(2r\eps)}{2r}\right)\right).
\]
The planar geodesic curvature of the $xy$-projection of $\Ga$ is
$h(\eps)=-2r$.  For any $\la\neq 0$ we consider the union of the
surfaces $\Sg_\la(\Ga)$ and $\widetilde{\Sg}_\la(\Ga)$ given in
Proposition~\ref{prop:sigmala} and Remark~\ref{re:reverse2},
respectively.  Easy computations from \eqref{eq:geocoor2} show that
the singular curves $\Ga_1$ and $\Ga_2$ are vertical translations of
$\Ga$ by $c_1(\la)\,T$ and $c_2(\la)\,T$, where
\begin{align*}
c_1(\la)&=\frac{s_\eps}{2\la}+
\frac{\text{sgn}(\la)\,\pi-2\la\,s_\eps}{4r^2}-\frac{(r^2+\la^2)\,
\sin(2\la\,s_\eps)}{4\la^2r^2},
\\
c_2(\la)&=\frac{\text{sgn}(\la)\,\pi}{2\la^2}-c_1(\la).
\end{align*}
In the first equation above $s_{\eps}$ is the cut function associated to
$\Ga_1$.  In general $\Ga_1\neq \Ga_2$ so that we can
extend the surface by geodesics of the same curvature orthogonal to
$\Ga_i$.  As indicated in Remark~\ref{re:reverse} and according with
the value of $\dot{\Ga}_i$, in order to preserve the constant mean
curvature $\la$ we must consider the surfaces
$\widetilde{\Sg}_{-\la}(\Ga_1)$ and $\Sg_{-\la}(\Ga_2)$.  Two new
singular curves $\Ga_{12}$ and $\Ga_{22}$ are obtained.  We repeat
this process by induction so that at any step $k+1$ we leave from the
singular curves $\Ga_{1k}$ and $\Ga_{2k}$ by the corresponding
orthogonal geodesics of curvature $(-1)^k\la$.  We denote by
$\mathcal{L}_\la$ the union of all these surfaces.  This is a $C^2$
immersed surface (in fact, it is $C^\infty$ outside the singular curves) and,
by construction, it is volume-preserving area-stationary with constant
mean curvature $\la$.  Any singular curve $\Ga_{ik}$ of
$\mathcal{L}_\la$ is a vertical translation of the helix $\Ga$ by the
vector $c_{ik}(\la)\,T$, where
\begin{align*}
c_{1k}(\la)&=k\,c_1(\la)-\text{sgn}(\la)\,\left[\frac{k}{2}\right]\,\frac{\pi}{2\la^2},
\\
c_{2k}(\la)&=\frac{\text{sgn}(\la)\,\pi}{2\la^2}-c_{1k}(\la),
\end{align*}
where $[x]$ denotes the greatest integer less than or equal to $x\in\rr$.

The singular curves $\Ga_{ik}$ could coincide depending on the values
of $\la$.  For example, an easy analytical argument shows that there
is a discrete set of values of $\la\in (0,r)$ for which $\Ga_1$
coincides with $\Ga$ (those for which $c_1(\la)$ is an integer
multiple of $\pi/(2r^2)$).  This situation is not possible when
$\la^2\geq r^2$.  In fact, for the case $r=\la=1$ explicit
computations from the equations above show that all the curves
$\Ga_{ik}$ are different.  So the resulting surface contains
infinitely many singular helices.  Also, it is not difficult to see
that for a discrete set of values of $\la\in (0,r)$, we have
$\Ga_{1i}=\Ga_{2i}$, so that we can obtain complete surfaces
$\mathcal{L}_\la$ with any given even number of singular curves.  In
general, the surfaces $\mathcal{L}_\la$ are not embedded.
\end{example}

In Theorem~\ref{th:classification} we will prove that any complete
volume-preserving area-stationa\-ry surface $\Sg$ in $\hh^1$ with
singular curves and non-vanishing mean curvature is congruent with one
of the surfaces $\mathcal{S}_\la$ or $\mathcal{L}_\la$ introduced
above.  We need the following strong restriction on the
singular curves of $\Sg$ obtained as a consequence of
Propositions~\ref{prop:c2singcurves} and \ref{prop:sigmala}~(iv).

\begin{theorem}
\label{th:curve}
Let $\Sg$ be a complete, oriented, $C^2$ immersed volume-preserving
area-stationa\-ry surface in $\hh^1$ with non-vanishing mean curvature.
Then any connected singular curve of $\Sg$ is a complete geodesic of
$\hh^1$.
\end{theorem}

\begin{proof}
Let $\mathcal{C}$ be a connected singular curve of $\Sg$.  By
Proposition \ref{prop:c2singcurves} we know that $\mathcal{C}$ is a
$C^2$ smooth horizontal curve.  We consider the unit normal $N$ to
$\Sg$ such that $N=T$ along $\mathcal{C}$.  Let $H$ be the mean
curvature of $\Sg$ with respect to $N$.  By using
Theorem~\ref{th:chmy} (ii) and Remark~\ref{rem:lambdaorientation}, for
any $p\in\mathcal{C}$ there is a small neighborhood of $p$ in $\Sg$
foliated by geodesics of curvature $H$ leaving from $\mathcal{C}$.  By
Theorem~\ref{th:constant} these geodesics are characteristic curves of
$\Sg$ and meet $\mathcal{C}$ orthogonally.

Let $\Ga$ be any closed arc of $\mathcal{C}$.  We parameterize $\Ga$
by arc-length $\eps\in [a,b]$.  By compactness we can find a small
$r>0$ such that, for any $\eps\in [a,b]$, the geodesic
$\ga_\eps:[0,r)\to\hh^1$ of curvature $H$ with initial conditions
$\Ga(\eps)$ and $J(\dot{\Ga}(\eps))$ is entirely contained in $\Sg$.
This implies that $\Sg$ and the surface $\Sg_{H}(\Ga)$ in
Proposition~\ref{prop:sigmala} locally coincides at one side of $\Ga$.
Moreover, as $\Sg$ is complete we deduce that $\Sg_H(\Ga)\sub\Sg$.  In
particular, $\Ga_1$ is a piece of a singular curve of $\Sg$ and so it
is $C^2$ smooth by Proposition \ref{prop:c2singcurves}.  As $\Sg$ is
volume-preserving area-stationary we deduce by Theorem
\ref{th:constant} that the geodesics $\ga_{\eps}$ meet $\Ga_1$
orthogonally.  This implies by Proposition \ref{prop:sigmala} (iv)
that the cut function $s_{\eps}$ is constant along $\Ga$.  As $\Ga$ is
an arbitrary closed arc of $\mathcal{C}$, we have proved that the
$xy$-projection of $\mathcal{C}=(x,y,t)$ is a straight line or a
planar circle.  Finally, by integrating the ``horizontal'' equation
$\dot{t}=\dot{x}y-x\dot{y}$ (as was done in
Section~\ref{sec:geodesics}) we conclude that $\mathcal{C}$ is a
complete geodesic of $\hh^1$.
\end{proof}

Now, we will see how to apply our previous results to describe all
compact vo\-lu\-me-preserving area-stationary surfaces in $\hh^1$.

The first relevant results about compact surfaces with constant mean
curvature in $\hh^1$ were given in \cite[Theorem~E]{chmy}, where it
was obtained an interesting restriction on the topology of an immersed
surface inside a spherical $3$-dimensional pseudo-hermitian manifold
under the weaker assumption that the mean curvature is bounded outside
the singular set.  The arguments in the proof use the local behavior
of the singular set studied in Theorem~\ref{th:chmy} and Hopf Index
Theorem for line fields.  They also apply to $\hh^1$ so that we get

\begin{proposition}[\cite{chmy}]
\label{prop:genus} Any compact, connected, $C^2$ immersed surface
$\Sg$ in $\hh^1$ with constant mean curvature is homeomorphic either
to a sphere or to a torus.
\end{proposition}

Moreover, in \cite[\S 7, Examples 1 and 2]{chmy} we can find examples
of constant mean curvature surfaces of spherical and toroidal type
inside the standard pseudo-hermitian $3$-sphere.  In $\hh^1$ we may
expect, by analogy with the Euclidean space, the existence of immersed
tori with constant mean curvature \cite{wente}.  However, this is not
possible as a consequence of our following result, that could be
interpreted as a counterpart in $\hh^1$ to Alexandrov uniqueness
theorem for embedded surfaces in $\rr^3$.

\begin{theorem}[Alexandrov Theorem in $\hh^1$]
\label{th:alexandrov}
Let $\Sg$ be a compact, connected, $C^2$ immersed volume-preserving
area-stationary surface in $\hh^1$.  Then $\Sg$ is congruent with a
sphere $\sph_{H}$ of the same constant mean curvature.
\end{theorem}

\begin{proof}
From the Minkowski formula \eqref{eq:minkowski} we have that the
constant mean curvature $H$ of $\Sg$ with respect to the inner normal
must be positive.  Observe that $\Sg$ must contain a singular point.
Otherwise Theorem \ref{th:ruled} would imply that $\Sg$ is foliated by
complete geodesics, a contradiction since any geodesic of $\hh^1$
leaves a compact set in finite time (Remark \ref{re:propertiesgeo}).
On the other hand $\Sg$ cannot contain a singular curve since this
curve would be a complete geodesic by Theorem \ref{th:curve} and $\Sg$
is compact.  We conclude that $\Sg$ has an isolated singularity.  We
finally invoke Theorem \ref{th:spheres} to deduce that $\Sg$ is
congruent with a sphere $\sph_{H}$ of the same mean curvature.
\end{proof}

Now, we shall prove the following classification theorem

\begin{theorem}
\label{th:classification}
Let $\Sg$ be a complete, oriented, connected, $C^2$ immersed
volume-preserving area-stationary surface in $\hh^1$ with
non-vanishing mean curvature.  If $\Sg$ contains a singular curve then
$\Sg$ is congruent either with the surface $\mathcal{S}_{H}$ in Example
\ref{ex:cilindros} or with the surface $\mathcal{L}_{H}$ in Example
\ref{ex:helices} of the same mean curvature as $\Sg$.
\end{theorem}

\begin{proof}
Let $\Ga$ be a connected horizontal curve of $\Sg$.  By Theorem
\ref{th:curve} we know that $\Ga$ is a complete geodesic of $\hh^1$.
After applying a left translation and a vertical rotation we can
suppose that $\Ga$ coincides either with the $x$-axis or with a helix
passing through the origin.  We can choose the unit normal $N$ to
$\Sg$ so that $N=T$ along $\Ga$.  By Theorem \ref{th:chmy} (ii) and
Remark~\ref{rem:lambdaorientation} there is $r>0$ such that the
geodesics $\ga_{\eps}:[0,r]\to\hh^1$ and
$\widetilde{\ga}_{\eps}:[0,r]\to\hh^1$ of curvature $H$ with initial
conditions $\Ga(\eps)$ and $J(\dot{\Ga}(\eps))$ (resp.  $\Ga(\eps)$
and $-J(\dot{\Ga}(\eps))$) are contained in $\Sg$.  As $\Sg$ is
complete and connected we can prolong these geodesics until they meet
a singular curve.  This implies that the union of the surfaces
$\Sg_{\la}(\Ga)$ and $\widetilde{\Sg}_{\la}(\Ga)$ constructed in
Proposition \ref{prop:sigmala} and Remark \ref{re:reverse2} is
included in $\Sg$.  The proof then follows by using the description of
the surfaces $\mathcal{S}_{\la}$ and $\mathcal{L}_{\la}$ in Examples
\ref{ex:cilindros} and \ref{ex:helices} together with the completeness
and the connectedness of $\Sg$.
\end{proof}

\begin{remark}
The previous result and Theorem \ref{th:spheres} provide the
description of complete $C^2$ immersed area-stationary surfaces under
a volume constraint in $\hh^1$ with non-empty singular set and
non-vanishing mean curvature.  Unduloids, cylinders and nodoids in
$\hh^1$ are examples of complete volume-preserving area-stationary
surfaces in $\hh^1$ with non-vanishing mean curvature and empty
singular set, see \cite{revolucion}.
\end{remark}

The arguments in this section can also be used to construct examples
and obtain restrictions on complete area-stationary surfaces in
$\hh^1$ with singular curves.

Let $\Ga=(x,y,t)$ be a $C^{k+1}$ ($k\geq 1$) horizontal curve in
$\hh^1$ parameterized by arc-length $\eps\in I$.  We denote by
$\ga_{\eps}:\rr\to\hh^1$ the geodesic of curvature zero and initial
conditions $\ga_\eps(0)=\Ga(\eps)$ and
$\dot{\ga}_\eps(0)=J(\dot{\Ga}(\eps))$.  We know from
Section~\ref{sec:geodesics} that $\ga_\eps$ is a horizontal straight
line.  We consider the map
$F(\eps,s)=\ga_{\eps}(s)=(x(\eps,s),y(\eps,s),t(\eps,s))$ given by
\begin{align}
\nonumber
x(\eps,s)&=x(\eps)-s\,\dot{y}(\eps),
\\
\label{eq:geocoor3}
y(\eps,s)&=y(\eps)+s\,\dot{x}(\eps),
\\
\nonumber
t(\eps,s)&=t(\eps)-s\,(x(\eps)\,\dot{x}(\eps)+y(\eps)\,\dot{y}(\eps)).
\end{align}
The Jacobi vector field $V_{\eps}(s):=(\ptl F/\ptl\eps)(\eps,s)$ along
$\ga_{\eps}$ can be computed from \eqref{eq:geocoor3} so that we get
\[
V_{\eps}(s)=(\dot{x}(\eps)-s\,\ddot{y}(\eps))\,X+(\dot{y}(\eps)+
s\,\ddot{x}(\eps))\,Y-s\,T.
\]
It follows that $\escpr{V_{\eps},T}<0$ on $(0,+\infty)$ and
$\escpr{V_\eps,T}>0$ on $(-\infty,0)$.  As a consequence the map
$F:I\times\rr\to\hh^1$ defines a complete immersed surface
$\Sg_{0}(\Ga)$.  By using Theorem~\ref{th:constant} we obtain that
$\Sg_0(\Ga)$ is a $C^k$ area-stationary surface whenever $k\geq 2$.
By following the proof of Theorem~\ref{th:classification} we deduce
the following geometric description of area-stationary surfaces
with singular curves.

\begin{proposition}
\label{prop:onesingu}
Let $\Sg$ be a complete, oriented, connected, $C^2$ immersed
area-stationary surface in $\hh^1$.  Then $\Sg$ contains at most one
singular curve $\Ga$.  In that case $\Sg$ consists of the union of
all the horizontal lines in $\hh^1$ orthogonal to $\Ga$.
\end{proposition}

The result above shows that the strong condition obtained in
Theorem~\ref{th:curve} does not hold for area-stationary surfaces.  We
can construct examples of such surfaces just by leaving from an
arbitrary horizontal curve by horizontal straight lines.  For example,
area-stationary helicoidal surfaces in $\hh^1$ are obtained when the
initial curve is a geodesic of non-zero curvature
\cite[Theorem~D]{pauls}.  Note that Proposition~\ref{prop:onesingu}
together with the already mentioned result in \cite{chmy} that any
complete minimal surface with an isolated singularity must coincide
with a Euclidean plane provides the complete description of complete
area-stationary surfaces in $\hh^1$ with non-empty singular set.

It is difficult to get a complete classification of minimal or
constant mean curvature surfaces without singular points in $\hh^1$,
see \cite{ch}.

We will say that a $C^1$ surface $\Sg$ is \emph{vertical} if the
vertical vector $T$ is contained in $T_{p}\Sg$ for any $p\in\Sg$.  A
complete vertical surface $\Sg$ is foliated by vertical straight
lines.  Since a $C^2$ vertical surface has no singular points, to have
constant mean curvature $H$ implies that $\Sg$ is either
area-stationary in case $H=0$, or volume-preserving area-stationary in
case $H\neq 0$.  From Theorem~\ref{th:ruled} is easy to get the
following, compare with \cite[Lemma~4.9]{gp},

\begin{proposition}
Let $\Sg$ be a $C^2$ complete, connected, immersed, oriented, constant
mean curvature surface in $\hh^1$.  If $\Sg$ is vertical then $\Sg$ is
either a vertical plane, or a right circular cylinder.
\end{proposition}

\section{\textbf{The isoperimetric problem in $\hh^1$}}
\label{sec:iso}
\setcounter{equation}{0}

The isoperimetric problem in $\hh^1$ consists of finding global
minimizers of the sub-Riemannian perimeter under a volume constraint.
For any Borel set $\Om\subeq\hh^1$ the \emph{perimeter} of
$\Om$ is defined by
\begin{equation*}
\pp(\Om):=\sup\,\left\{\int_\Om\divv(U)\,dv;\, |U|\leq 1\right\},
\end{equation*}
where the supremum is taken over $C^1$ \emph{horizontal vector fields}
with compact support on $\hh^1$.  In the definition above, $dv$ and
$\divv(\cdot)$ are the Riemannian volume and divergence of the left
invariant metric $g$, respectively.  This notion of perimeter coincides
with the $\hh^1$-perimeter introduced in \cite{cng} and \cite{fsc}.
For a set $\Om$ bounded by a surface $\Sg$ of class $C^2$ we have
$\pp(\Om)=A(\Sg)$ by virtue of the Riemannian divergence theorem.

It is not difficult to prove that the perimeter is $3$-homogeneous
with respect to the family of dilations in \eqref{eq:dilations}, see
for instance \cite[Lemma 4.5]{msc}.  Precisely, for any Borel set
$\Om\subeq\hh^1$ and any $s\in\rr$ we have
\[
V(\varphi_s(\Om))=e^{4s}\,V(\Om),\qquad
\pp(\varphi_s(\Om))=e^{3s}\,\pp(\Om).
\]
This property leads us to the isoperimetric inequality
\begin{equation}
\label{eq:isoineq}
\pp(\Om)^4\geq\alpha\,V(\Om)^3,
\end{equation}
that holds for any Borel set $\Om\subeq\hh^1$.  Inequality
\eqref{eq:isoineq} was first obtained by P.~Pansu \cite{pansu2} for
regular sets.  Many other generalizations have been established but
always without the sharp constant $\alpha$, see \cite{garonhieu} and 
\cite{dgn2}.

An \emph{isoperimetric region} in $\hh^1$ is a set $\Om\sub\hh^1$ such that
\[
\pp(\Om)\leq\pp(\Om')
\]
amongst all sets $\Om'\sub\hh^1$ with $V(\Om)=V(\Om')$.

The existence of isoperimetric regions was proved by G. P. Leonardi
and S. Rigot \cite[Theorem~2.5]{lr} in the more general context of
Carnot groups, see also \cite[Theorem~13.7]{dgn}.  We summarize their
results in the following theorem.

\begin{theorem}[\cite{lr}]
\label{th:leorig} For any $V>0$ there is an isoperimetric region $\Om$
in $\hh^1$ with $V(\Om)=V$.  The set $\Om$ is, up to a set of measure
zero, a bounded connected open set.  Moreover, the boundary $\ptl\Om$
is Alhfors regular and verifies condition $B$.
\end{theorem}

The condition B in the theorem above is a certain separation property.
It means that there is a constant $\beta>0$ such that for any
Carnot-Carath\'eodory ball $B$ centered on $\ptl\Om$ with radius
$r\leq 1$ there exist two balls $B_1$ and $B_2$ with radius $\beta r$
such that $B_1\sub B\cap\Om$ and $B_2\sub B-\overline{\Om}$.

The properties in Theorem~\ref{th:leorig} are not sufficient to
describe the isoperimetric regions in $\hh^1$.  In 1983 P.~Pansu made
the following

\vspace{0,2cm} \noindent\textbf{Conjecture} (\cite[p.~172]{pansu3}).
In the Heisenberg group $\hh^1$ any isoperimetric region bounded by
a smooth surface is congruent with a sphere $\sph_\la$.

In the last years many authors have tried to adapt to the Heisenberg
setting different proofs of the classical isoperimetric inequality in
Euclidean space.  In \cite{monti2}, \cite{monti} and \cite{leomas} it
was shown that there is no a direct counterpart in $\hh^1$ to the
Brunn-Minkowski inequality in Euclidean space, with the consequence
that the Carnot-Carath\'edory metric balls in $\hh^1$, cannot be the
solutions.  Recently, expecting that symmetrization could work in
$\hh^1$, interest has focused on solving the isoperimetric
problem restricted to certain sets with additional symmetries.  It has
been recently proved by D. Danielli, N. Garofalo and D.-M. Nhieu that
the sets $\Om_\la$ bounded by the spherical surfaces $\sph_\la$ are
the unique solutions in the class of sets bounded by two $C^1$ graphs
over the $xy$-plane \cite[Theorem~1.1]{dgn2}.  An intrinsic
description of the solutions was given by G. P. Leonardi and S. Masnou
\cite[Theorem~3.3]{leomas}, where it was proved that any sphere
$\sph_\la$ is the union of all the geodesics of curvature $\la$ in
$\hh^1$ connecting the poles.  In \cite{revolucion} we pointed out
that assuming $C^2$ smoothness and rotationally symmetry of
isoperimetric regions, these must be congruent with the spheres
$\sph_\la$.  We also mention the interesting recent work \cite{bc} in
which it is proved that the flow by mean curvature of a $C^2$ convex
surface in $\hh^1$ described as the union of the radial graphs
$t=\pm f(|z|)$, with $f'> 0$, converges to the spheres $\sph_{\la}$.

\vspace{0,1cm} The regularity of isoperimetric regions in $\hh^1$ is
still an open question.  The regularity of the spheres $\sph_\la$ and
of the examples of complete volume-preserving area-stationary surfaces
in Section~\ref{sec:mainresult} may suggest that the isoperimetric
solutions in $\hh^1$ are $C^\infty$ away from the singular set and
only $C^2$ around the singularities.

By assuming $C^2$ regularity of the solutions we can use the
uniqueness of spheres in Theorem~\ref{th:alexandrov} to solve the
isoperimetric problem in $\hh^1$.

\begin{theorem}
\label{th:iso}
If $\Om$ is an isoperimetric region in $\hh^1$ bounded by a $C^2$
smooth surface $\Sg$, then $\Om$ is congruent with a set bounded by a
sphere $\sph_\la$.
\end{theorem}

\begin{proof}
Let $\Om$ be an isoperimetric region of class $C^2$ in $\hh^1$.  By
using Theorem~\ref{th:leorig} we can assume that $\Om$ is bounded and
connected.  The boundary $\Sg=\ptl\Om$ is a $C^2$ compact surface with
finitely many connected components.  Let us see that $\Sg$ is
connected.  Otherwise we may find a bounded component $\Om_{0}$ of
$\hh^1-\overline{\Om}$.  Consider the set $\Om_{1}=\Om\cup\Om_0$.  It
is clear that $V(\Om_1)>V(\Om)$ and $\pp(\Om_1)<\pp(\Om)$.  Thus by
applying an appropriated dilation to $\Om_{1}$ we would obtain a new
set $\Om'$ so that $V(\Om')=V(\Om)$ and $\pp(\Om')<\pp(\Om)$, a
contradiction since $\Om$ is isoperimetric.  As $\Sg$ is a $C^2$
compact, connected, volume-preserving area-stationary surface in
$\hh^1$, we conclude by Alexandrov (Theorem~\ref{th:alexandrov}) that
$\Sg$ is congruent with a sphere $\sph_\la$.
\end{proof}

\begin{remark}[The isoperimetric constant in $\hh^1$]
The area of the sphere $\sph_\la$ can be easily computed from
\eqref{eq:spheregraphs}. Using polar coordinates and
Fubini's theorem we get
\[
A(\sph_\la)=\frac{\pi^2}{\la^3}.
\]
On the other hand, we can use Minkowski formula~\eqref{eq:minkowski}
to compute the volume of the set $\Om_\la$ enclosed by $\sph_\la$.  We
obtain
\[
V(\Om_\la)=\frac{3\pi^2}{8\la^4}.
\]
In case the $C^2$ regularity of isoperimetric sets in $\hh^1$ was
established, we would deduce from Theorem~\ref{th:iso} that the
optimal isoperimetric constant in \eqref{eq:isoineq} would be given by
\[
\alpha=\frac{A(\sph_\la)^4}{V(\Om_\la)^3}=\bigg(\frac{8}{3}\bigg)^3\pi^2.
\]
\end{remark}

\providecommand{\bysame}{\leavevmode\hbox to3em{\hrulefill}\thinspace}
\providecommand{\MR}{\relax\ifhmode\unskip\space\fi MR }
\providecommand{\MRhref}[2]{%
   \href{http://www.ams.org/mathscinet-getitem?mr=#1}{#2} }
   \providecommand{\href}[2]{#2}

\end{document}